\documentclass[12pt]{amsart}
\usepackage{amsmath, amscd, amssymb, amsthm,amsfonts}
\usepackage[nohug]{diagrams}
\usepackage[T1]{fontenc}
\makeatletter

\usepackage[all]{xy}
\usepackage{dbnsymb}
\usepackage{stmaryrd}

\ifx\pdfoutput\undefined
\usepackage{graphicx}
\else
\usepackage[pdftex]{graphicx}
\usepackage{epstopdf}
\fi

\DeclareMathOperator{\In}{In}
\DeclareMathOperator{\Out}{Out}

\DeclareMathOperator{\Tout}{Tout}
\DeclareMathOperator{\Tin}{Tin}

\DeclareMathOperator{\BBar}{Bar}
\DeclareMathOperator{\Cobar}{Cobar}

\DeclareMathOperator{\Har}{Harrison}
\DeclareMathOperator{\Torus}{Torus}
\DeclareMathOperator{\pTorus}{pre-Torus}
\DeclareMathOperator{\Gr}{Gr}

\DeclareMathOperator{\Val}{val}
\DeclareMathOperator{\sign}{sgn}

\DeclareMathOperator{\Lie}{L}
\DeclareMathOperator{\Ass}{A}
\DeclareMathOperator{\Comm}{C}

\DeclareMathOperator{\Group}{Group}
\DeclareMathOperator{\colim}{colim}

\DeclareMathOperator{\tree}{Tree}

\DeclareMathOperator{\Endo}{End}

\DeclareMathOperator{\Diff}{Diff}

\DeclareMathOperator{\Top}{Top}
\DeclareMathOperator{\Cell}{Cell}

\DeclareMathOperator{\Chain}{Kom}

\DeclareMathOperator{\Bij}{Bij}
\DeclareMathOperator{\Morph}{Hom}

\DeclareMathOperator{\Aut}{Aut}
\DeclareMathOperator{\Hpty}{Htpy}

\DeclareMathOperator{\Sh}{Sh}
\DeclareMathOperator{\Obj}{Ob}

\newcommand{\inp}[1]{\ensuremath{\langle #1 \rangle}}
\newcommand{\inpp}[1]{\ensuremath{\langle\langle #1 \rangle\rangle}}
\newcommand{\module}[0]{\operatorname{-mod}}

\newcommand{\cal}[1]{\ensuremath{\mathcal{#1}}}

\newcommand{\normaltext}[1]{\textnormal{#1}}

\newcommand{\CPic}[1]{
\begin{minipage}{.45in}
\includegraphics[scale=.75]{#1}
\end{minipage}
}

\newcommand{\MPic}[1]{
\begin{minipage}{.35in}
\includegraphics[scale=.45]{#1}
\end{minipage}
}

\newcommand{\BPic}[1]{
\begin{minipage}{1in}
\includegraphics[scale=1.5]{#1}
\end{minipage}
}

\theoremstyle{plain}
\newtheorem{theorem}[subsection]{Theorem}
\newtheorem*{observation}{Observation}

\newtheorem{corollary}[subsection]{Corollary}
\newtheorem{lemma}[subsection]{Lemma}
\theoremstyle{remark}
\newtheorem{remark}[subsection]{Remark}

\theoremstyle{definition}

\newtheorem{definition}[subsection]{Definition}

\newtheorem*{notation}{Notation}

\numberwithin{equation}{section}

\newcommand{\Xp}[0]{\ensuremath{X_{v}}}
\newcommand{\Xr}[0]{\ensuremath{X}}

\newcommand{\Xpq}[0]{\ensuremath{BH_{v}}}
\newcommand{\Yp}[0]{\ensuremath{|NL_{v}|}}
\newcommand{\Yq}[0]{\ensuremath{|NL_{w}|}}
\newcommand{\Ypq}[0]{\ensuremath{|NL_{v\#w}|}}
\newcommand{\Lp}[0]{\ensuremath{NL_{v}}}

\begin{document}
\title{Topological Field Theories and Harrison Homology}

\author[Benjamin Cooper]{Benjamin Cooper}

\address{Department of Mathematics, University of Virginia, Charlottesville, VA 22904}
\email{bjc4n\char 64 virginia.edu}

\begin{abstract}
Tools and arguments developed by Kevin Costello are adapted to families of
``Outer Spaces'' or spaces of graphs. This allows us to prove a version of
Deligne's conjecture: the Harrison homology associated to a homotopy
commutative algebra is naturally a module over a cobordism category of
3-manifolds.
\end{abstract}

\maketitle

\section{Introduction}\label{intro}
%\footnote{MSC Codes: 81T45, 18D50, 57M50}
%\footnote{Keywords: homology operations, operads, C infinity algebras, outer
%  space}

A recent theorem of Kevin Costello (\cite{MR2298823}) illustrates the
relationship between homotopy associative algebras ($\Ass_{\infty}$
algebras) and moduli spaces of Riemann surfaces. See also
\cite{KSNotes1}. It is an exciting addition to a story which began when
Deligne conjectured that the action of the homology of configuration spaces
on the Hochschild homology of an associative algebra, $HH_*(A,A)$, lifts to
an action defined at the chain level. Deligne's conjecture was shown to be
true (see \cite{MR1898414}), but thinking of configuration spaces as moduli
of genus 0 surfaces leads to a more general theorem: the chain level action
of genus 0 surfaces extends to a natural action by surfaces of all
genera. More specifically, the chain complex computing the Hochschild
homology of an $\Ass_{\infty}$ algebra is the object associated to the
circle by a 2-dimensional topological field theory.

Costello considers the moduli spaces of open, open-closed and closed Riemann
surfaces, and defines three categories related by inclusions,
$$j : \cal{O} \hookrightarrow \cal{OC} \hookleftarrow \cal{C} : i .$$

The categories of modules over the open, open-closed and closed categories
are open, open-closed and closed topological field theories.  The category
of open theories can be identified with the category of cyclic
$\Ass_{\infty}$ algebras. Given such an algebra $A$, the inclusions $i$ and
$j$ yield a functor
$$i^* \circ j_* : \cal{O}\module \to \cal{C}\module,$$

which assigns to $A$ the closed topological field theory $i^*\circ
j_*(A)$. This closed theory associates a chain complex, $i^*\circ
j_*(A)(S^1)$, to the circle. The homology of this chain complex is
isomorphic to Hochschild homology of $A$, $H_*(i^*\circ j_*(A)(S^1)) \cong
HH_*(A,A)$. It follows that the geometrically defined maps in $\cal{C}$
determine natural operations on the Hochschild complex of $A$.

%It is then determined by direct computation that the object associated to a
%circle must compute Hochschild homology, $H_*(i^*\circ j_*(A)(S^1)) \cong
%HH_*(A,A)$.

Recent work by Hatcher, Vogtmann and Wahl (\cite{MR2113904,MR2174267})
suggests that natural choices of open, open-closed and closed categories may
be obtained from the classifying spaces of mapping class groups of doubled
handlebodies, or 3-manifolds of the form $\#^gS^1\times S^2\#^e D^3\#^t S^1
\times D^2$. Such classifying spaces can be modelled by spaces of metric
graphs. In this paper, the open, open-closed and closed categories are
categories of rational chains on these spaces.

We first prove that the category of modules over the open category is
equivalent to the category of cyclic $\Comm_{\infty}$ algebras. The
extension from the open category to the open-closed category yields a
Costello-type theorem:

\begin{theorem}
  There exists a category $\mathcal{OC}$ of 3-manifolds with objects given
  by boundary spheres $S$ and tori $T$. There are open and closed
  subcategories with inclusions,
$$j : \mathcal{O} \hookrightarrow \mathcal{OC} \hookleftarrow \cal{C} : i.$$

The category of cyclic $\Comm_\infty$ algebras is equivalent to the category
of modules over the category $\mathcal{O}$. For any $\Comm_\infty$ algebra
$A$, the homology of the object associated to the torus is isomorphic to the
Harrison homology of $A$.
$$H_*(i^*\circ j_*(A))(T) \cong \Har_*(A,A)$$

%In particular, the chain complex computing the Harrison homology of $A$ is
%naturally a module over the closed category $\cal{C}$.
In particular, the closed category acts naturally on a chain complex
computing the Harrison homology of $A$.
\end{theorem}

In section \ref{extandtorussec}, the chain complex $\Torus(A) = (i^*\circ
j_*)(A)(T)$ will be described explicitly. From this description, it will
follow that the homology, $H_*(\Torus(A))$, agrees with the Harrison homology
$\Har_*(A,A)$ of $A$.

The spaces of graphs appearing in this paper are natural extensions of the
``Outer Spaces'' which originally appeared in \cite{MR830040}. Connections
between Outer Space and the homotopy commutative operad appear in
\cite{MR1341841, MR1301191, GetzKapMod, MR2431663}. The construction in this
paper can be viewed as a version of Deligne's conjecture in the ``classical
limit'' corresponding to the homotopy commutative operad in Kontsevich's
``three worlds'' \cite{MR1247289}.

The author would like to thank Justin Roberts for his encouragement,
Nathalie Wahl for her helpful emails during the first writing of this
document, Jim Stasheff for his interest and the referees for their helpful
comments and suggestions.

\section{Algebra and Operads} \label{algebra and operads}

The underlying ring in all constructions will be the field of rational
numbers. We denote by $\Top$ the category of topological spaces, by $\Group$
the category of groups and by $\Chain$ the category of chain complexes.

\subsection{Monoidal Categories} \label{monoidal categories}

A category $\mathcal{C}$ is \emph{symmetric monoidal} if it is equipped with
a bifunctor,
$$-\otimes- : \mathcal{C} \times \mathcal{C} \to \mathcal{C},$$

an object $1$ and isomorphisms,

\begin{enumerate}
\item $(a \otimes b) \otimes c \cong a \otimes (b \otimes c) $
\item $ 1 \otimes a \cong a \cong a \otimes 1$
\item $ a \otimes b \cong b \otimes a $
\end{enumerate}

satisfying coherence conditions, see \cite{MR1712872}.  There are monoidal
structures on $\Top$, $\Group$ and $\Chain$ given by disjoint union, product
and tensor product respectively.

A \emph{monoidal} functor $F : \mathcal{C} \to \mathcal{D}$ between
symmetric monoidal categories is equipped with maps $F(a)\otimes F(b) \to
F(a\otimes b)$ that are natural in both $a$ and $b$, and satisfy
associativity and commutativity criteria.

\begin{definition}{($\Obj(\cal{C})$)} \label{objcat}
Every symmetric monoidal category \cal{C} has a subcategory $\Obj(\cal{C})$
with the same objects. The morphisms of $\Obj(\cal{C})$ are generated by permutations of tensors.
$$a \otimes a' \cong a' \otimes a$$
\end{definition}

Notice that $\Obj(\Obj(\cal{C})) = \Obj(\cal{C})$. This notation agrees with
\cite{MR2298823}.

\subsection{Differential Graded Categories} \label{dg categories}

All of the categories in this paper will have extra structure in a sense
that can be captured by the idea of enrichment. A category $\mathcal{C}$ is
\emph{enriched} over a monoidal category $\mathcal{D}$ when, for all $X$,
$Y \in \Obj(\mathcal{C})$,
$$\Morph_{\mathcal{C}}(X,Y) \in \Obj(\mathcal{D})$$

and the composition in $\cal{C}$ respects this $\cal{D}$ structure:
\begin{align*}
\circ : \Morph_{\cal{C}}(a,b) \otimes \Morph_{\cal{C}}(b,c) &\to \Morph_{\cal{C}}(a,c)\\ & \in \Morph_{\cal{D}}(\Morph_{\cal{C}}(a,b) \otimes \Morph_{\cal{C}}(b,c),\Morph_{\cal{C}}(a,c)),
\end{align*}

for all $ a,b,c\in \Obj(\cal{C})$. A category $\cal{C}$ that satisfies
$\Morph_{\cal{C}}(X,Y)\in \Top$ will be called a \emph{topological
  category}.  For instance, the sets $\Morph_{\Top}(X,Y)$ can be endowed
with the compact open topology.
%%  showing that $\Top$ is enriched over $\Top$.  A category in which
%% $\Morph_{\cal{C}}(X,Y)\in \Top$ will be called \emph{topological}.
%% A \emph{linear} category is a category enriched over
%%$\Vect$. 

If $\cal{C}$ is enriched over $\cal{D}$ and $F : \cal{D} \to \cal{E}$ is a
monoidal functor, then there is a category $F_*(\cal{C})$ enriched over
$\cal{E}$ defined by:
$$\Obj(F_*(\cal{C})) = \Obj(\cal{C}) \quad\normaltext{ and }\quad F_*(\cal{C})(a,b) = F(\cal{C}(a,b)).$$

For example, the functor $B : \Group \to \Top$ giving the classifying space
of a group determines the functor $B_*$ which maps categories enriched over
$\Group$ to categories enriched over $\Top$. Another important example is
$C_*(-;\mathbb{Q})$; rational singular chains. If $\mathcal{C}$ is a
topological category then there is a category $C_*(\mathcal{C};\mathbb{Q})$
defined by:
\begin{align*}
\Obj(C_*(\mathcal{C};\mathbb{Q})) &= \Obj(\mathcal{C})\\ 
\Morph_{C_*(\mathcal{C};\mathbb{Q})}(A,B) &=C_*(\Morph_{\mathcal{C}}(A,B);\mathbb{Q}).
\end{align*}

%If our category is enriched in cellular spaces then the functor
A version of this functor will be defined for cellular spaces in section
\ref{natchains}. Categories of the form $C_*(\cal{C};\mathbb{Q})$ are
examples of differential graded categories.

A \emph{differential graded} or \emph{dg} category is a category enriched
over $\Chain$.  A \emph{differential graded symmetric monoidal}, or
\emph{dgsm} category, is a symmetric monoidal category which is differential
graded. The category $\Chain$ is an example of a dgsm category. If
$\mathcal{C}$ is a dg category then $H_*(\mathcal{C})$ is a category
enriched over the category of graded vector spaces.

A dgsm functor or \emph{morphism} of dgsm categories, $F: \mathcal{C} \to
\mathcal{D}$, is a monoidal functor which respects the differential graded
structure. This is a monoidal functor of categories enriched over $\Chain$,
as defined above.  Two dgsm categories $\cal{C}$ and $\cal{D}$ are
\emph{quasi-isomorphic} when there is a dgsm functor $F : \cal{C} \to \cal{D}$
such that $H_*(F)$ is full, faithful and induces isomorphisms on objects.

\subsection{Modules Over Differential Graded Categories}\label{dg modules}

If $\mathcal{C}$ is a dgsm category then a \emph{left $\mathcal{C}\module$}
is a dgsm functor $\mathcal{C} \to \Chain$. A \emph{right
  $\mathcal{C}\module$} is a dgsm functor $\mathcal{C}^{op} \to
\Chain$. Note that as functors, modules must respect the differential
graded structure, specifically if
$$F_{a,b} : \Morph_{\cal{C}}(a,b) \to \Morph_{\Chain}(F(a),F(b))$$

then $d\circ F_{a,b} = F_{a,b} \circ d$ for all $a,b \in \Obj(\cal{C})$.

Maps between modules $M$ and $N$ are natural transformations $\phi : M \to
N$ of the underlying functors which satisfy the following conditions:

\begin{enumerate}
\item All of the elements $\phi(a) \in \Morph(M(a),N(a))$ are chain maps.
\item The natural transformation $\phi$ respects the monoidal structure,
\begin{center}
$$\begin{diagram}
  M(a)\otimes M(a') & \rTo & N(a)\otimes N(a')\\
\dTo        &                   & \dTo\\
M(a \otimes a') & \rTo & N(a \otimes a').\\
\end{diagram}$$
\end{center}
\end{enumerate}

The category of left (right) modules over $\mathcal{C}$ will be denoted by
$\mathcal{C}$-mod (mod-$\mathcal{C}$). 

For a functor to be monoidal we \emph{only} require the existence of a map
$$F(a) \otimes F(b) \to F(a\otimes b)$$

satisfying the axioms described in section 2.1. It is often the case that
these maps satisfy stronger conditions. A module is said to be \emph{split}
when the monoidal structure maps $F(a) \otimes F(b) \to F(a \otimes b)$ are
isomorphisms and \emph{h-split}, or homologically split, when they are
quasi-isomorphisms. For instance, a TQFT in the sense of Atiyah \emph{is} a
split module over the cobordism category.
%% \cite{MR1001453}.

The usual product of categories extends to one which respects the dgsm
structure. If $\mathcal{C}$ and $\mathcal{D}$ are categories then there is a
category $\mathcal{C} \otimes \mathcal{D}$ defined by 
\begin{align*}
\Obj(\mathcal{C} \otimes \mathcal{D}) &= \Obj(\cal{C}) \times \Obj(\cal{D})\\
\Morph_{\mathcal{C} \otimes \mathcal{D}}(a\times c, b\times d) &= \Morph_{\mathcal{C}}(a,c) \otimes \Morph_{\mathcal{D}}(b,d).
\end{align*}

If $\mathcal{C}$ and $\cal{D}$ are differential graded then
$\cal{C}\otimes\cal{D}$ is differential graded using the usual tensor
product of chain complexes. If $\cal{C}$ and $\cal{D}$ are monoidal then
$\cal{C}\otimes\cal{D}$ is monoidal using $(a\times c)\otimes (b\times d) =
(a\otimes b)\times (c\otimes d)$.

If $\mathcal{C}$ and $\mathcal{D}$ are dgsm categories then a
$\mathcal{D}-\mathcal{C}$ \emph{bimodule} is a dgsm functor from the
category $\mathcal{D} \otimes \mathcal{C}^{op}$ to $\Chain$. 

The following observation will be used to define the bimodule $\mathcal{OC}$
appearing in theorem \ref{octoocg} section \ref{hvwtheoremsec}.

\begin{observation}\label{canonicalmoduleobs}
Every dgsm category $\mathcal{C}$ yields a $\mathcal{C}-\mathcal{C}$
bimodule, $\cal{C} : \cal{C}\otimes\cal{C}^{op} \to \Chain$ given by
$\mathcal{C}(x\times y) = \Morph_{\mathcal{C}}(y,x) $.
\end{observation}

If $M$ is a $\mathcal{D}-\mathcal{C}$ bimodule and $N$ is a left
$\mathcal{C}\module$ then there exists a left $\mathcal{D}\module$, $M
\otimes_\mathcal{C} N$, defined on objects $b\in \Obj(\mathcal{D})$ by,
$$(M \otimes_\mathcal{C} N)(b) = \bigoplus_{a\in\Obj(\mathcal{C})} M(b,a)\otimes N(a),$$

modulo the relation, $\sim$, which makes the diagram below to commute,
\begin{center}$$\begin{diagram}
M(b,a)\otimes \Morph_{\mathcal{C}}(a',a) \otimes N(a') & \rTo & M(b,a) \otimes N(a)\\
\dTo        &                   & \dTo\\
M(b,a') \otimes N(a') & \rTo & (M \otimes_{\mathcal{C}} N)(b).\\
\end{diagram}$$\end{center}

Explicitly,
$$f^*(g) \otimes h \sim g \otimes f_*(h),$$

for $f\in\Morph_{\cal{C}}(a',a), g\in M(b,a) \normaltext{ and } h\in
N(a')$.

Although dgsm modules do not form a dg category, they do possess a reasonable
notion of weak equivalence. A map $\varphi : M \to M'$ between $M,
M'\in\mathcal{C}$-mod is a \emph{quasi-isomorphism} if $\varphi_* : H_*(M(a))\to
H_*(M'(a))$ for all $a \in \Obj(\cal{C})$.

A functor $F : \mathcal{C} \to \mathcal{D}$ between categories of modules is
\emph{exact} when it preserves quasi-isomorphisms.  Two functors $F,G :
\mathcal{C} \to \mathcal{D}$ are \emph{quasi-isomorphic}, $F \simeq G$, when
there are natural transformations $\varphi : F \to G$ such that $\varphi(c)$
is a quasi-isomorphism for all $c \in \Obj(\mathcal{C})$. Two categories
$\mathcal{C}$ and $\mathcal{D}$ are \emph{isomorphic} or
\emph{quasi-equivalent}, $\mathcal{C} \cong \mathcal{D}$ if there are
functors $F : \mathcal{C} \to \mathcal{D}$ and $G : \mathcal{D} \to
\mathcal{C}$ such that $FG \simeq 1$ and $GF \simeq 1$.

A module $M$ is \emph{flat} when the functor $-\otimes M$ is exact.  Since
most of the constructions to follow will involve considering dgsm categories
and their modules up to quasi-isomorphism, strictly speaking, we should be
working in a derived category. As such the tensor product $M\otimes N$ of a
$\cal{D}-\cal{C}$ bimodule $M$ and a left $\cal{C}\module$ $N$ as above
should be defined by $M\otimes^{\mathbb{L}}_{\cal{C}} N = M
\otimes_{\cal{C}} \BBar_{\cal{C}} N$ where $\BBar_{\cal{C}} N$ is $N$
tensored with the $\BBar$ construction on $\cal{C}$. This gives a canonical
flat replacement (see \cite{MR2298823}).

Any dgsm functor $F : \cal{C} \to \cal{D}$ between dgsm categories
determines a pair of functors between the corresponding categories of
modules, $F^* : \cal{D}\module \to \cal{C}\module$ and $F_* : \cal{C}\module
\to \cal{D}\module$. The functor $F^*$ is called restriction and $F_*(M) =
\cal{D}\otimes_{\cal{C}} M$ is called the induction functor. The latter is
defined using the tensor product above and the $\cal{D}-\cal{C}$ bimodule
structure on $\cal{D}$ inherited from $F$.

\begin{theorem}{(\cite{MR2298823})}\label{modeq} If $F : \mathcal{C} \to \mathcal{D}$ is a quasi-isomorphism of dgsm categories, then the induction and restriction functors,
$$\mathbb{L}F_* : \cal{C}\module \rightleftarrows \cal{D}\module : F^*$$

are inverse quasi-isomorphisms between the categories of left (right)
$\cal{C}$ modules and left (right) $\cal{D}$ modules respectively.
\end{theorem}

\subsubsection{Cellular Chains}\label{natchains}

If $X$ is a cellular space then we would like the equivalence
$C_*^{cell}(X;\mathbb{Q})\simeq C_*(X;\mathbb{Q})$ to be natural.  In order
to accomplish this, our chain complexes are defined to be a colimit over all
maps from cellular spaces into a given space (see \cite{MR2298823}).

A \emph{cellular space} $X$ is a CW complex of finite type; in other words,
there are finitely many cells in each dimension.  In particular, each cell
attaches to only finitely many other cells. If $X^i$ is the $i$-skeleton of
$X$ then $f : X\hookrightarrow Y$ is a map of cellular spaces when it is
continuous and $f^{-1}(Y^i) = X^i$. Let $\Cell \subset \Top$ be the
subcategory of cellular spaces and cellular maps. For any topological space
$Y$, define
$$C_*(Y;\mathbb{Q}) = \underset{X\in\Cell\downarrow Y}{\colim} C_*^{cell}(X;\mathbb{Q})$$

where $\Cell\downarrow Y$ is the over category and
$C^{cell}_*(-;\mathbb{Q})$ denotes the functor given by taking rational
cellular chains. It follows that if $Y$ is a cellular space then the map
$C^{cell}_*(Y;\mathbb{Q})\to C_*(Y;\mathbb{Q})$ is natural.

\subsection{Operads}\label{operadssubsec}

After a brief discussion of operads and cyclic operads, we introduce the
$\BBar$ and $\Cobar$ functors and define the associative and
associative commutative operads: $\Ass$ and $\Comm$. The operad
$\Comm_{\infty}$ will first be introduced as a quotient of the
$\Ass_{\infty}$ operad. In section \ref{barcobar}, $\Comm_{\infty}$ will be
defined in terms of the $\Cobar\circ\BBar$ construction.

\subsubsection{Operads} \label{operadsubsubsec}
In what follows operads will be used to encode axioms for various kinds of
algebras and to control stratifications of certain spaces of graphs. For
more information regarding operads, see \cite{MR1436914,
  MR1898414,MR2131012}.

A differential graded \emph{operad} $\mathcal{O}$ is a sequence of chain
complexes $\{\mathcal{O}(n)\}_{n=1}^{\infty}$ and composition maps
$$ \gamma : \mathcal{O}(k) \otimes \mathcal{O}(n_1)\otimes\cdots\otimes\mathcal{O}(n_k) \to \mathcal{O}(n_1+\cdots+n_k)$$

together with an action of the symmetric group $\Sigma_{n}$ on
$\mathcal{O}(n)$ and a unit $1\in \mathcal{O}(1)$. The compositions $\gamma$
are required to be $\Sigma$-equivariant, associative and unital. In all
cases to follow, chain complexes will be finite dimensional and
$\mathcal{O}(1)$ will be one dimensional.

A map of operads $\varphi : \mathcal{O} \to \mathcal{O}'$ is given by a
collection $\{\varphi_n : \cal{O}(n)\to \cal{O'}(n)\}_{n=1}^{\infty}$ of
$\Sigma$-equivariant chain maps that commute with the operad compositions
and take units to units.  Two operads $\cal{O}$ and $\cal{O'}$ are
\emph{quasi-isomorphic} when there is a map $\varphi : \cal{O} \to
\cal{O'}$, the individual components of which induce isomorphisms on
homology.

Cooperads are operads with the arrows, $\gamma$, reversed. There is a
completely analogous category of differential graded cooperads, see \cite{GetzJones}.

Given a chain complex $X$, define the \emph{endomorphism operad},
$\Endo_X$, by
$$\mathcal{\Endo}_X(n) = \Morph^*_{\Chain}(X^{\otimes n}, X).$$

Composition is given by composition of chain maps and the action of $\Sigma_n$
is given by permuting the arguments of $f \in \Endo_X(n)$. A chain complex
$X$ is an \emph{algebra} over an operad $\mathcal{O}$ when there is a morphism
of operads $\mathcal{O}\to \Endo_X$.

A differential graded \emph{cyclic operad} is an operad $\cal{O} =
\{\mathcal{O}(n)\}_{n=1}^{\infty}$ such that the action of $\Sigma_n$ on
$\cal{O}(n)$ lifts to an action of $\Sigma_{n+1}$ on $\cal{O}(n)$. An
algebra $X$ over a cyclic operad $\cal{O}$ is required to possess a
non-degenerate bilinear form which is invariant with respect to the
operations of $\cal{O}$, see \cite{MR1358617}. 

Operads are usually pictured as rooted trees with vertices labelled by some
distinguished collection of symbols. The composition $\gamma$ corresponds to
gluing the roots of $k$ such trees to the unrooted edges of a single tree
with $k+1$ boundary edges.  A cyclic operad is an operad in which the trees
representing operations lack a preferred root. Cyclic operations can be
manipulated in the plane, see section \ref{rel to dg alg}.

\subsubsection{Homotopy operads} \label{examples of operads}

In this section we give explicit models for the operads $\Comm$, $\Ass$,
$\Comm_{\infty}$ and $\Ass_{\infty}$.  The usual definition of
$\Comm_{\infty}$ is given as a quotient of $\Ass_{\infty}$ by the shuffle
relations. Since dg operads defined by quotients cannot control moduli
spaces, such as those found in section \ref{outerspaces}, in section
\ref{barcobar}, the $\Cobar\circ\BBar$ functor is introduced in order to
remove the shuffle relations.

The \emph{commutative operad} $\Comm = \{\Comm(n)\}_{n=1}^{\infty}$ is both
the main object of interest and the simplest operad:
$$\Comm(n) = \mathbb{Q} \quad\normaltext{ for all }\quad n \geq 1,$$

concentrated in degree 0. If $X$ is a vector space then $X$ is an algebra
over the commutative operad when $X$ is an associative commutative
algebra. $\Comm$ extends to a cyclic operad. A cyclic $\Comm$ algebra is an
associative commutative algebra $X$ with an inner product $\inp{-,-} :
X\otimes X \to \mathbb{Q}$ that satisfies,
$$\inp{a\cdot b,c} = \inp{a,b\cdot c}.$$

In other words, $X$ is a commutative Frobenius algebra.

The $\Ass_{\infty}$ operad is generated by all possible compositions of
$n$-fold operations $m_n$ subject to the relation that
$$ \partial m_n(1,\ldots,n) = \sum_{\substack{i+j = n+1\\i,j \geq 2}} \sum_{s=0}^{n-j} (-1)^{j+s(j+1)} m_i(1,\ldots,m_j(s+1,\ldots,s+j+1),\ldots,n),$$

where $m_n(1,\ldots,n)$ is the operation $m_n$ labelled by its $n$
inputs. The degree of $m_n$ is $n-2$.

Elements of the operad $\Ass_{\infty}$ are usually pictured as rooted trees
in the plane in which the $n+1$-valent vertices represent the operation
$m_n$. The operation $m_n$ is sometimes represented by a disk with $n$
distinct boundary points and an extra boundary point corresponding to the
root.  In this case, a composition of the form
$m_i(1,\ldots,m_j(\ldots),\ldots,n)$ is represented by two such disks glued
together along one of their boundary points.

The homotopy associative commutative or $\Comm_{\infty}$ operad is usually
introduced as a quotient of the $\Ass_{\infty}$ operad by shuffle relations.
The operad $\Comm_{\infty}$ is the kernel of the map
$\Ass_{\infty}\to \Lie_{\infty}$ obtained by extending the map $A \to L$
defined by $[a,b] = ab - ba$.

A $(p,q)$-\emph{shuffle}, $\sigma \in \Sh(p,q)$, is a permutation $\sigma \in
\Sigma_{p+q}$ which satisfies,
$$\begin{array}{ccc}
\sigma(1) < \sigma(2) < \ldots < \sigma(p) &  \normaltext{ and }   & \sigma(p+1) < \sigma(p+2) < \ldots < \sigma(p+q).
\end{array}$$

The $\Comm_{\infty}$ operad is obtained from the $\Ass_{\infty}$ operad by imposing the relations,
$$ \sum_{\sigma \in \Sh(i,n-i)} \sign(\sigma) m_n({\sigma(1)},\ldots,{\sigma(n)}) = 0 ,$$

for all $1 < i < n$, where $\sign(\sigma)$ is the sign of a permutation. For
instance, when $k=2$, the relation becomes,
$$ m_2(a, a') - m_2(a', a) = 0. $$

Cyclic $\Comm_\infty$ and $\Ass_\infty$ algebras possess a non-degenerate
inner product $\inp{-,-}$ which satisfies
$$ \inp{m_n(x_0,\ldots,x_{n-1}),x_n} = (-1)^{(n+1) |x_0| \sum_{i=1}^{n-1} |x_i|}\inp{m_n(x_1,\ldots,x_n),x_0} .$$

If $M$ is an $n$-manifold then the de Rham complex $\Omega^*(M)$ is an
example of a $\Comm_{\infty}$ algebra. The transfer theorem determines a
$\Comm_{\infty}$ algebra structure on $H^*(M;\mathbb{R})$. If $M$ is also
compact then $H^*(M;\mathbb{R})$ is cyclic; the inner product is the duality
pairing, see \cite{MR1103672} and \cite{MR2503530}.

There is a map of operads $\alpha : \Comm_{\infty} \to \Comm$, defined by
$$\alpha(m_2) = m_2 \quad\normaltext{ and }\quad \alpha(m_j) = 0 \normaltext{ if } j \ne 2$$

which is a quasi-isomorphism. We'd like to think of $\Comm_{\infty}$ as a
free resolution of $\Comm$. Unfortunately, since we have added the shuffle
relations, $\Comm_{\infty}$ is not free in the appropriate sense. In order
to obtain a dg operad homotopy equivalent to $\Comm$, which is free of
relations, we introduce the $\Cobar\circ\BBar$ functor in section
\ref{barcobar}.

\subsection{Resolutions of operads} \label{resolutions of operads}

In this section we introduce definitions for graphs and use these
definitions to construct the $\BBar$ and $\Cobar$ functors. The Bar
construction is a functor which takes a dg operad $\cal{P}$ to a dg cooperad
$\BBar(\cal{P})$, while the Cobar construction is a functor taking a dg
cooperad $\cal{O}$ to a dg operad $\Cobar(\cal{O})$. These form an
adjunction between the categories of operads and cooperads, the unit of
which,
$$\eta_{\cal{O}} : \cal{O} \to \Cobar(\BBar(\cal{O}))$$

is a quasi-isomorphism of operads. 

\subsubsection{Graphs}\label{gorprelim}

A \emph{graph} $G$ is a finite set that has been partitioned in two
ways: into pairs $e = \{a,b\}$ called \emph{edges} and into sets $H(v) =
\{h_1,h_2,\ldots,h_n\}$ called \emph{vertices}.
$$ G = \coprod_{e} \{a,b\} = \coprod_{v} H(v)$$

Denote the set of vertices of $G$ by $V(G)$ and the set of edges of $G$ by
$E(G)$. The elements of $G$ will be called \emph{half edges}. Two half edges
$a,b\in G$ \emph{meet} if $a,b\in H(v)$ for some vertex $v$. Given an edge
$e \in E(G)$, the set $e = \{x,y\}$ is the set of half edges associated to
$e$ in $G$. For each vertex $v\in V(G)$, the set $H(v)$ is the set of half
edges associated to $v$ in $G$. The \emph{valence} $\Val(v)$ of $v\in V(G)$
is the number of half edges or $|H(v)|$. All graphs in this document are
required to have vertices $v$ of valence $\Val(v) = 1 $ or $\Val(v) \geq 3$
unless otherwise noted.

Two graphs $G$ and $H$ are \emph{isomorphic} when there is a bijective set map
between half edges $\varphi : H \to G$ that respects the two partitions.

A \emph{subgraph} $H$ of $G$ is a graph formed by the set of all vertices of
$G$ together with some subset of the set of edges of $G$.  A \emph{cycle} of
$G$ is a subgraph $C\subset G$ given by an ordered sequence of edges which
begin and end at the same vertex.

The \emph{boundary} $\partial(G)$ of a graph $G$ is the collection of edges
that contain a vertex having valence one. An \emph{internal} edge is an edge
not in the boundary while an \emph{external} edge is not internal. 

Let $[n]$ be the set $\{1,\ldots,n\}$. A graph $G$ is \emph{boundary
  labelled} if there is a choice of partition $\partial(G) = \In(G) \cup
\Out(G)$ of the boundary into a set of \emph{incoming} and \emph{outgoing}
edges together with bijections $i_G : [|\In(G)|] \to \In(G)$ and $o_G :
[|\Out(G)|] \to \Out(G)$. A notion of boundary labelling for geometric
graphs appears in section \ref{hty eq grp}.

A \emph{geometric graph} is a 1-dimensional CW complex. Every graph $G$ has
an associated geometric graph, $|G|$, in which the $0$-skeleton is given by
the vertices $V(G)$ and the $1$-skeleton is formed by gluing on $1$-cells
corresponding to the edges. We may refer to graphs as either combinatorial
or geometric when it is necessary to make a distinction.

% It will be clear whether by $|X|$ we mean the
%cardinality of the set $X$ or the CW complex associated to a graph $X$.

A graph $G$ is \emph{connected} when $H_0(|G|) \cong \mathbb{Q}$. A graph $G$
has \emph{genus} $g$ if $H_1(|G|) \cong \mathbb{Q}^g$. A \emph{forest} is a
graph of genus $0$. A \emph{tree} is a connected forest. A \emph{rooted
  tree} is a tree together with a choice of outgoing edge, the rest of the
boundary edges being incoming. A tree with a single internal vertex will be
called a \emph{corolla}. An $n$-$\tree$ is a tree with $n$ incoming edges.

Given an edge $e \in E(G)$, $e = \{x,y\}$, we can form a new graph $G/e$ by
removing $e$ and replacing $H(x)$ and $H(y)$ with $H(x)\cup H(y) -
\{x,y\}$. This operation, called \emph{edge collapse}, is a homotopy
equivalence of $|G|$ if $x$ and $y$ are not contained in the same set of
half edges $H(v)$.  Collapsing a forest $F\subset G$ is called \emph{forest
  collapse}.

\subsubsection{Orientations}\label{orprelim}

If $V_*$ is a graded vector space then the \emph{$j$-fold desuspension}
$V[j]_*$ is given by $V[j]_i = V_{i+j}$. An \emph{orientation} of a graded
vector space $W$ of dimension $n = \dim(W)$ is a non-zero vector in the
exterior algebra $\det(W) = \Lambda^n(W)[-n]$. The dual is defined by
$\det(W)^* = \Lambda^n(W)[n]$. If $S$ is a set then we orient $S$ using
$\det(S) = \det(\mathbb{Q}\inp{S})$. Two orientations are equivalent when they
are positive scalar multiples of each other. An \emph{orientation of a
  graph} $G$ is defined to be an element of
$$\det(G) = \det(E(G)) \otimes \det(\Out(G))\otimes \det(H_0(G)) \otimes \det(H_1(G))^*[O-\chi]$$

where $O$ is the number of outgoing edges and $\chi = \chi(G)$ is the Euler
characteristic of $G$. Using this convention, a graph is placed in degree
$|E(G)|$.  There are maps,
$$\det(G_0)\otimes \det(G_1) \to \det(G_0\#G_1)\quad\normaltext{and}\quad \det(G_0\coprod G_1) \cong \det(G_0) \otimes \det(G_1).$$

\subsubsection{The $\BBar$ and $\Cobar$ constructions}\label{barcobar}

If $S$ is a set and $\mathcal{O}$ is a cyclic dg (co)operad then a
\emph{labelling} of $S$ by $\mathcal{O}$ is defined by the coinvariants
trick:
$$\mathcal{O}(S) = (\mathcal{O}(n) \times \Bij([n+1],S))_{\Sigma_{n+1}}$$

where $\Bij([n+1],S)$ is the set of bijections from $S$ to $[n+1]=
\{1,\ldots,n+1\}$ and $\Sigma_{n+1}$ acts diagonally. If $T$ is a tree then a \emph{labelling} of $T$ by $\mathcal{O}$ is determined by assigning to each vertex
$v$ an element of $\mathcal{O}(H(v))$,
$$\mathcal{O}(T) = \bigotimes_{v\in V(T)} \mathcal{O}(H(v)) .$$

The collapse of an internal edge $c : T \to T/e$ induces maps of
labellings. If we denote by $e$ the vertex obtained by the edge collapse and
by $v$ and $w$ the two identified end points then there are maps,
$$\mathcal{O}(\Val(v)) \otimes \mathcal{O}(\Val(w)) \to \mathcal{O}(\Val(e)) \quad\normaltext{ and }\quad \mathcal{P}(\Val(e)) \to \mathcal{P}(\Val(v)) \otimes \mathcal{P}(\Val(w)).$$

(Recall from \ref{gorprelim} that $\Val(v)$ is the valence of the vertex
$v$.)  Tensoring the above with identity yields maps $c_* :
\mathcal{O}(T)\to \mathcal{O}(T/e)$ and $c^* : \mathcal{P}(T/e)\to
\mathcal{P}(T)$. These maps $c_*$ and $c^*$ are used to define the $\BBar$
and $\Cobar$ differentials below.

The \emph{Bar construction} $\BBar(\mathcal{O})$ of a cyclic differential
graded operad $\mathcal{O}$ is the dg cooperad of labelled trees with an
edge contracting differential. Explicitly,

$$ \BBar(\mathcal{O})(n) = \bigoplus_{\substack{n-\tree T\\|T|=1}} \mathcal{O}(T)\otimes\det(T) \gets \cdots \gets \bigoplus_{\substack{n-\tree T\\|T|=n-1}} \mathcal{O}(T)\otimes\det(T) .$$

The \emph{Cobar construction} $\Cobar(\mathcal{P})$ of a cyclic differential
graded cooperad $\mathcal{P}$ is the dg operad of labelled trees with an
edge expanding differential. Concretely,

$$ \Cobar(\mathcal{P})(n) = \bigoplus_{\substack{n-\tree T\\|T|=1}} \mathcal{P}(T)\otimes\det(T)^* \to \cdots \to \bigoplus_{\substack{n-\tree T\\|T|=n-1}} \mathcal{P}(T)\otimes\det(T)^* .$$

In the formulas above $|T|$ is the number of internal vertices of $T$. The
complex is graded so that the term spanned by trees with one internal vertex
is situated in degree $0$. Alternatively, the grading is determined by the
orientation, see section \ref{orprelim}.

The differential $\delta$ either contracts or expands edges. It can be
described by its matrix elements, $(\delta)_{T,T'}$.  If $T'$ is not
isomorphic to $T/e$ for some internal edge $e\in T$ then the corresponding component
of $\delta$ is set to zero. Otherwise, let $c : T \to T' \cong T/e$ so that
if $c_* : \mathcal{O}(T) \to \mathcal{O}(T')$ or $c^* : \mathcal{P}(T') \to
\mathcal{P}(T)$ are the maps above then $\delta$ is given by
$$ (\delta)_{T,T'} = c_* \otimes p_e \quad\normaltext{ or }\quad (\delta)_{T',T} = c^* \otimes p^e .$$

%%$$(\delta)_{T,T'} : \mathcal{O}(T)\otimes\det(T) \to \mathcal{O}(T')\otimes\det(T')\normaltext{ and } (\delta)_{T',T} : \mathcal{P}(T')\otimes\det(T')^*\to\mathcal{P}(T)\otimes\det(T)^*.$$
%%$$\mathcal{O}(T)\otimes\det(T) \xrightarrow{(\delta)_{T,T'}} \mathcal{O}(T')\otimes\det(T')\normaltext{  or  } \mathcal{P}(T')\otimes\det(T')^*\xrightarrow{(\delta)_{T',T}}\mathcal{P}(T)\otimes\det(T)^*.$$

If collapsing the edge $e$ identifies the vertices $u$ and $v$ to a vertex
$e$, then the map of orientations $p_e : \det(T) \to \det(T')$ is given by,
$$p_e(y_0\wedge \cdots \wedge e \wedge \cdots\wedge y_n) = y_0 \wedge \cdots\wedge \hat{e}\wedge \cdots \wedge y_n$$

and the orientation map $p^e$, coupled with the expanding differential is
defined analogously. In either case, if the operad $\cal{O}$, or cooperad
$\cal{P}$, has a non-trivial differential then the total differential is the
sum of the differential defined above together with the original
differential.

The composition for the operad $\Cobar(\cal{P})$ is given by grafting
boundary edges and eliminating the resulting bivalent vertex. This satisfies
the Leibniz rule with respect to the differential defined above. Notice that
$\Cobar(\mathcal{P})$ is generated by $\mathcal{P}$-labelled corolla.

The $\BBar$ and $\Cobar$ functors form an adjunction. The counit and unit maps of this adjunction,
$$\BBar(\Cobar(\cal{P}))\to \cal{P} \quad\normaltext{ and }\quad \cal{O} \to \Cobar(\BBar(\cal{O})),$$
are quasi-isomorphisms, see \cite{GetzJones}.

\subsection{Relation to Differential Graded Algebra} \label{rel to dg alg}

The language of differential graded operads and their algebras in section
\ref{operadssubsec} is an important special case of the language of
differential graded categories and their modules, see section \ref{dg
  categories}.  In this section we establish a connection between sections
\ref{operadssubsec} and \ref{dg categories}.

%On one hand, operads appear as nice presentations for dg
%categories obtained from topological categories of moduli spaces; on the
%other hand, for every a dg operad \cal{O} there is a dgsm category
%$\cal{O}^{\flat}$ so that the category of h-split $\cal{O}^{\flat}$ modules
%is isomorphic to the category of $\cal{O}$ algebras.

Given a dg operad $\mathcal{O}$, we can define the \emph{enveloping category}
$\mathcal{O^{\flat}}$ to be the dgsm category generated by one object $X$
and morphisms generated by
$$ \Morph_{\cal{O}^{\flat}}(X^{\otimes n},X) = \cal{O}(n)$$

using the monoidal structure. Pictorially, if operations $x\in \cal{O}(k)$
are represented by trees then $y \in \Morph_{\cal{O}^{\flat}}(X^{\otimes n},
X^{\otimes m})$ is a disjoint union of trees. By construction, the category
$\cal{O}^{\flat}$ includes factorization isomorphisms,
$$\theta_{n,m} = \Morph_{\cal{O}^{\flat}}(X^{\otimes n}, X^{\otimes m}) \cong \bigotimes_{i=1}^m \Morph_{\cal{O}^{\flat}}(X^{\otimes n_i}, X)\quad\textnormal{such that}\quad \sum_i n_i = n .$$

Maps of operads induce functors between their associated enveloping
categories. The following is an immediate consequence of the above
construction.

\begin{lemma}
  The category of $\mathcal{O}$-algebras is equivalent to the category of
  split left $\mathcal{O^{\flat}}$ modules.
\end{lemma}
\begin{proof}
  Any functor $F: \mathcal{O^{\flat}} \to \Chain$ identifies the object $X$
  with a chain complex $F(X)$ and, by split monoidality, identifies the object
  $X^{\otimes m}$ with $F(X)^{\otimes m}$. Consider the action of
  $\mathcal{O^{\flat}}$ on $F(X)$. Using the factorization map $\theta_{n,m}$, 
$$\varphi =  \varphi_1\otimes \cdots\otimes\varphi_m\quad\normaltext{ such that }\quad \varphi_i : X^{\otimes n_i}\to X$$ 
where $n = n_1 + n_2 + \cdots + n_m$. Each map $\varphi_i$ is also an element of
$\cal{O}(n_i)$. This identification commutes with the categorical composition of $\cal{O}^\flat$ and the operadic composition of $\cal{O}$.
\end{proof}

Split modules do not behave as well under quasi-isomorphism as h-split
modules. The next lemma tells us that, for our purposes, these two notions
of split are equivalent.

\begin{lemma}
  There is an equivalence of categories between the category of h-split left
  $\cal{O}^{\flat}$ modules and the category of split left
  $\cal{O}^{\flat}$ modules.
\end{lemma}
\begin{proof}
  The equivalence is determined by a functor $\eta$ from h-split to split
  modules. If $F$ is an h-split $\cal{O}^{\flat}$ module then define a split
  module $\eta(F)$ by $\eta(F)(X^{\otimes n}) = F(X)^{\otimes n}$.

Since $F$ is h-split there are quasi-isomorphisms $\varphi_{X^j}:
\eta(F)(X^{\otimes j}) \to F(X^{\otimes j})$. We must extend $\eta(F)$ to a
functor. Each $m_j \in \cal{O}(j)$ induces a map, $(m_j)_* : F(X)^{\otimes
  j} \to F(X)$. These are natural with respect to the maps
$\varphi_{X^j}$. For any $f \in \Morph_{\cal{O}^{\flat}}(X^{\otimes
  m},X^{\otimes n})$, $f = \theta^{-1}_{n,m}(m_{n_1}\otimes\cdots\otimes
m_{n_k})$.

So the action of $\cal{O}$ can be extended to an action of
$\cal{O}^{\flat}$, giving a unique split $\cal{O}^{\flat}$ module,
$\eta(F)$, which is quasi-isomorphic to the h-split $\cal{O}^{\flat}$ module
$F$ via $\{\varphi\}$.
\end{proof}

The following lemma allows us to simplify some rather complicated looking
operads. See section \ref{dg modules} for the definition of quasi-equivalent.

\begin{lemma} \label{envelopingequiv}
  If $\cal{O}_1$ and $\cal{O}_2$ are quasi-isomorphic operads then the
  associated enveloping categories $\cal{O}_1^{\flat}$ and
  $\cal{O}_2^{\flat}$ are quasi-isomorphic.
$$\cal{O}^{\flat}_1 \cong \cal{O}^{\flat}_2$$

In particular, it follows that the associated categories of modules are
quasi-equivalent.
$$\cal{O}_1^{\flat}\module \cong \cal{O}_2^{\flat}\module$$
\end{lemma}

The statement about modules follows from the lemmas and Theorem \ref{modeq}.

Cyclic differential graded operads $\mathcal{O}$ also yield dgsm categories
$\mathcal{O}^{\flat}$ with one object $X$ and morphisms generated by 
$$ \Morph_{\mathcal{O}^{\flat}}(X^{\otimes n},X) = \mathcal{O}(n), $$

together with cap and cup morphisms corresponding to an invariant inner product and its dual,
$$\inp{-,-} \in \Morph_{\mathcal{O}^\flat}(X\otimes X,\mathbb{Q})\quad\textnormal{and}\quad\inp{-,-}^* \in \Morph_{\mathcal{O}^\flat}(\mathbb{Q},X\otimes X).$$

These are represented by pictures,

$$\MPic{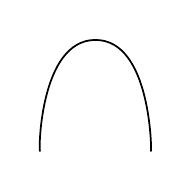}\quad\quad\normaltext{and}\quad\quad\MPic{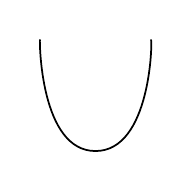},$$

which are subject to the S-bend relations:

$$\CPic{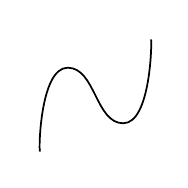}\quad = \quad \CPic{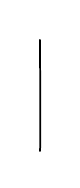} = \CPic{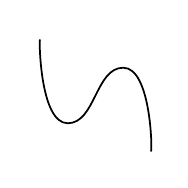}\quad.$$

The addition of caps and cups yields much larger morphism spaces;
$\Morph_{\mathcal{O}^{\flat}}(X^{\otimes n},X^{\otimes m})$ is now a space
of graphs (not a space of trees). Analogues of the previous lemmas
hold for $\cal{O}^{\flat}$ after $\cal{O}$ algebras are replaced by cyclic
$\cal{O}$ algebras.

\begin{remark}\label{reltomodularremark}
A differential graded PROP is a symmetric monoidal category which is
generated by a single object $x$ and enriched in the category of chain
complexes, see \cite{MR1898414, MR1712872}. The construction $-^\flat$ is a
functor from the category of cyclic dg operads to the category of dg PROPs.

Each dg modular operad $\mathcal{M}$ (see \cite{GetzKapMod}) determines a dg PROP $P\mathcal{M}$ where
$$\Morph_{P\mathcal{M}}(x^{\otimes n},x^{\otimes m}) = \oplus_g \mathcal{M}(g, n+m)$$
and the composition is constructed by gluing the corresponding collections of end points using the structure maps,
$$\circ_{ij} : \mathcal{M}(g,m)\otimes \mathcal{M}(g',n) \rightarrow \mathcal{M}(g+g',m+n-2).$$

A cyclic dg operad $\mathcal{O}$ determines a modular operad
$\mathcal{MO}$. Some authors refer to $\mathcal{MO}$
as the na\"{i}ve modular closure of $\mathcal{O}$, see \cite{refloch}
section 3.2.  If $\mathcal{O}$ is a cyclic dg operad then the PROP
$\mathcal{O}^\flat$ agrees with $P\mathcal{MO}$,
$$\mathcal{O}^\flat \cong P\mathcal{MO}.$$
\end{remark}

\section{$3$-dimensional Cobordism Categories} \label{3dcob}

In this section we define a dgsm category $\mathcal{M}$ called the
\emph{differential graded cobordism} category. A 3-dimensional topological
field theory will be a left \cal{M} module. In section \ref{ocandoc}, we
define the open, open-closed and closed subcategories of $\cal{M}$ which
will be used throughout the remainder of the paper. Although our focus is on
3-manifolds, the categories $\cal{N}$ and $\cal{M}$ have analogues in any
dimension.

Suppose that $M$ is a smooth manifold and let $\Diff(M,\partial)$ be the
group of diffeomorphisms of $M$ which fix a regular neighborhood of the
boundary. The \emph{mapping class group} $\Gamma(M,\partial)$ of $M$ is
defined to be $\pi_0\Diff(M,\partial)$. 

A \emph{labelled} surface $X$ is a surface $X$ together with a bijection
$[m] \to \pi_0(X)$ where $m = |\pi_0(X)|$ and $[m] = \{1,\ldots,m\}$. A
3-manifold with labelled boundary is analogous to a boundary labelled graph,
see sections \ref{gorprelim} and \ref{hty eq grp}.

%% Let $M$ be a smooth manifold with boundary. $\Diff(M)$ the group of
%% orientation preserving diffeomorphisms of $M$. Let $\Diff(M,\partial)
%% \subset \Diff(M)$ be the subgroup of diffeomorphisms which fix a regular
%% neighborhood of the boundary $\partial M$. The \emph{mapping class group}
%% $\Gamma(M,\partial)$ of $M$ is defined to be $\pi_0\Diff(M,\partial)$.

\begin{definition}{(\cal{N})}
The \emph{cobordism category} is a topological category $\cal{N}$ with
objects given by disjoint unions of orientable labelled surfaces.

A morphism $M' \in \Morph_{\mathcal{N}}(X,Y)$ is a triple $M' = (M, i, j)$
where $M$ is a diffeomorphism class (rel $\partial$) of smooth oriented
3-manifold whose boundary $\partial M = I \coprod J$ splits into a disjoint
union of labelled incoming surfaces $I$ and labelled outgoing surfaces $J$,
the orientations of which are induced by that of $M$. The maps $i : N(I) \to
X \times [0,\epsilon)$ and $j : N(J) \to Y \times [0,\epsilon)$ parameterize
    regular neighborhoods, $N(I), N(J) \subset M$, of the boundary. Any two
    choices of $\epsilon > 0$ define equivalent categories. For related
    discussion, see \cite{Stong}.

Given $A' = (A,i,j) \in \Morph_{\cal{N}}(X,Y)$ and $B' = (B, l, m) \in
\Morph_{\cal{N}}(Y,Z)$ define $B'\circ A' \in \Morph(X,Z)$ by
gluing: if $A\#B = A \coprod B / (x \sim y \normaltext{ if } j(x) = l(y))$
then $B'\circ A' = (A\#B,i,m)$. Associativity follows from the local nature of the
gluing composition.  Identity morphisms are given by thickened surfaces, $Y\times [0,1]$.

The category $\cal{N}$ has a symmetric monoidal structure given by disjoint
union. 
%%, which are completely fixed by the relevant diffeomorphisms.
\end{definition}

\begin{definition}{(\cal{M})}
The \emph{differential graded cobordism category} $\mathcal{M}$ is the
category of singular chains on classifying spaces of mapping class
groups of morphisms in $\cal{N}$. Specifically,
$$\Obj(\mathcal{M}) = \Obj(\mathcal{N}) \quad\normaltext{ and }\quad\Morph_{\mathcal{M}}(X,Y) =C_*(B\Gamma(\Morph_{\mathcal{N}}(X,Y),\partial);\mathbb{Q}).$$

We apply these functors to the triplets above in the most straightforward
way. If $M' = (M,i,j)$ is a morphism in $\mathcal{N}$ then
$\Gamma(M',\partial) = (\Gamma(M,\partial),i,j)$ and gluing of triples in
$\mathcal{N}$ as defined above induces a composition.

Specifically, if $A' = (A,i,j) \in \Morph_{\mathcal{N}}(X,Y)$, $B' = (B,l,m)
\in \Morph_{\mathcal{N}}(Y,Z)$ then given $(\phi,i,j) \in
\Gamma(A',\partial)$ and $(\psi,l,m) \in \Gamma(B',\partial)$, by requiring
that group elements fix a neighborhood of the boundary it follows that there
exists a map $\psi \# \phi : A\#B \to A\#B$ induced by $(\psi,\phi) :
A\coprod B \to A \coprod B$ so that $(\psi \# \phi,i,m)$ is a morphism in
$\Morph_{\Gamma(\mathcal{N},\partial)}(X,Z)$. The local nature of the gluing
implies associativity of the composition.

If $\Gamma(M',\partial) \in \Morph_{\Gamma(\mathcal{N},\partial)}(X,Y)$
then we say that $g \in \Gamma(M',\partial) = (\Gamma(M,\partial),i,j)$ when
$g \in \Gamma(M,\partial)$. Such elements form a group and so the functor
$B$ can be applied to $\Morph_{\Gamma(\mathcal{N},\partial)}(X,Y)$. We
apply $C_*(-;\mathbb{Q})$ to these classifying spaces. As discussed in
section \ref{dg categories}, both $B$ and $C_*(-;\mathbb{Q})$ are monoidal.
\end{definition}

Notice that the category $\mathcal{N}$ can be recovered as
$H_0(\mathcal{M};\mathbb{Q})$. We may think of $\mathcal{M}$ as a choice of
chain level representative for $\mathcal{N}$. Better terminology might be
level 0 differential graded cobordisms.
% because taking mapping class groups removes higher homotopical
%information from $B\Diff(M)$.

\begin{definition}{(TFT)}\label{tft}
A \emph{3-dimensional topological field theory} is an h-split left
$\mathcal{M}$ module. 
\end{definition}

\subsection{Open, Open-Closed and Closed Subcategories}\label{ocandoc}

The category $\mathcal{M}$ appears to be a very complicated object. We will
leverage the relationship between several much simpler subcategories of
$\mathcal{M}$: the \emph{open} category $\mathcal{O}$, the
\emph{open-closed} category $\mathcal{OC}$, and the \emph{closed} category
$\mathcal{C}$.

Let $\mathcal{S}$ be a collection of compact oriented 3-manifolds with
boundary. If $\inp{\mathcal{S}}$ is the subcategory of $\cal{N}$ generated
by $\mathcal{S}$ then a subcategory $\inpp{\mathcal{S}}$ of $\mathcal{M}$ is
\emph{generated by} $\mathcal{S}$ when $\inpp{\mathcal{S}}$ is
$C_*(B\Gamma(\inp{\mathcal{S}},\partial);\mathbb{Q})$.
%and
%$\inp{\mathcal{S}}$ is the subcategory of $\cal{N}$ generated by
%$\mathcal{S}$.

The categories below will use \emph{doubled handle bodies} with sphere and
torus boundary as generating manifolds. Let,
$$M_{(g,e,t)} = \#^gS^1\times S^2\#^e D^3\#^t S^1 \times D^2$$

be the connected sum of $g$ copies of $S^1\times S^2$, $e$ copies of $D^3$
and $t$ copies of $S^1\times D^2$. Notice that each $D^3$ summand introduces
a boundary 2-sphere and each $S^1\times D^2$ introduces a boundary
torus. The boundary of $M_{(g,e,t)}$ consists of $e$ 2-sphere and $t$
tori.
% Alternatively, $M_{(g,e,t)}$ can be constructed by starting with a
%doubled handle body of genus $g$, $\#^g S^1\times S^2$ then removing $e$
%solid three balls and $t$ solid tori.

We will adopt the following vector subscript notation for the remainder of
the paper.

\begin{notation}{$(M_v)$}\label{subscriptnote}
We write $M_v$ where $v = (g,i+j,n+m)$ for a manifold $M$ of
genus $g$ with labelled boundary consisting of $i$ incoming spheres, $n$
incoming tori, $j$ outgoing spheres and $m$ outgoing tori. An operation $\#$
is defined on composable subscripts by gluing the outgoing boundary of $M_v$
to the incoming boundary of $M_w$.
$$M_{v\#w} \cong M_v\#M_w$$
%If $w = (g',j+l,m+r)$ then
%% where $v\#w = (g+g'+j-1+m,i+l,n+r)$ since gluing the $j$ spheres together
%% adds $j-1$ factors of $S^1\times S^2$ and gluing the $m$ tori together adds
%% an additional $m$ factors of $S^1\times S^2$. 
\end{notation}

\begin{definition}{(\cal{OC})}
The \emph{open-closed category} $\mathcal{OC}$ is the subcategory
$\inpp{\mathcal{S}} \subset \mathcal{M}$ generated by $\mathcal{S} = \{
M_{(g,e,t)} \}$ such that there is always incoming and outgoing boundary. If
$t = 0$ then $e\geq 2$ and if $e=0$ then $t\geq 2$. In particular, when $t
\ne 0$, we require that there is always an incoming torus. The set
$\mathcal{S}$ is closed under composition.
\end{definition}

The open and closed categories are subcategories of the open-closed category.

\begin{definition}{(\cal{O} and \cal{C})}\label{opcatdef}
The \emph{open category} $\mathcal{O}$ is defined to be the subcategory of
$\mathcal{OC}$ whose objects are spheres and whose morphisms are generated by
the spaces $M_{v}$ where $v = (g,i+j,0)$. Similarly, the \emph{closed
  category} $\mathcal{C}$ is the subcategory of $\mathcal{OC}$ whose objects
are tori and whose morphisms are generated by the spaces $M_{v}$ where $v =
(g,0,n+m)$, (note $n\geq 1$).
%( \cite{MR2298823}).
\end{definition}

In each case, the composition is induced from gluing along boundaries
and identity morphisms are added as above.

\begin{definition}{(open, open-closed, closed TFT)}\label{octft}
 An \emph{open-closed} topological field theory is an h-split left
 $\mathcal{OC}$ module. An \emph{open} topological field theory is an
 h-split left $\mathcal{O}$ module. A \emph{closed} topological field theory
 is an h-split left $\mathcal{C}$ module.
\end{definition}

\section{Outer Spaces} \label{outerspaces}

In this section we will use the work of Hatcher, Vogtmann and Wahl on spaces
of graphs to reduce the categories $\cal{O}$ and $\cal{OC}$ to combinatorial
objects. In section \ref{hvwtheoremsec}, we show that mapping class groups
of the doubled handlebodies $M_v$ appearing in section \ref{ocandoc} are
rationally equivalent to certain groups associated to graphs. In section
\ref{opouterspacesec}, we construct ``Outer Spaces'' (see \cite{MR830040,
  MR2077676}) which model the rational homotopy type of the classifying
spaces of these groups. The associated group homology has been studied by
Hatcher and Vogtmann (\cite{MR1671188}) and is computed by the forested
graph complex. In section \ref{foresttocobar}, we show that this complex is
generated by a version of the $\Comm_{\infty}$ operad.

With the idea of ``classical degeneration'' in mind, it might be more
natural to consider the cobordism category of abstract tropical curves
\cite{MR2480499,MR2404949}. What follows is mostly speculation. An extension
of tropicalization to families of curves would allow one to define a
tropical analogue of conformal field theory and a restriction functor from
conformal field theories to tropical conformal field theories. Working with
chain complexes would remove the strict dependence of such theories on the
underlying geometry.  Work on tropical analogues of Gromov-Witten
theories, enumerative tropical geometry, suggests the existence of
topological tropical field theories. The references above suggest that
aspects of the material developed in this paper may coincide with such a
theory. The relationship between cyclic $\Comm_{\infty}$ algebras and the
rational homotopy theory of manifolds could allow one to compare tropical
degenerations of topological conformal field theories with constructions in
manifold theory.

\subsection{Homotopy Equivalence Groups}\label{hty eq grp}

We will now use a construction of Hatcher and Wahl \cite{MR2174267} to show
that the mapping class group of morphisms in the open, open-closed and
closed categories can be identified with automorphism groups of graphs.

A \emph{boundary torus} or \emph{balloon} is the geometric graph formed from
two edges with both ends of one edge glued to one end of the other. Define
the graph $G_v$ to be the geometric graph consisting of a wedge of $g$
circles with $e$ edges and $t$ boundary tori glued to the one base vertex
along the ends of their free edges.

\begin{center}
 \includegraphics[scale=.45]{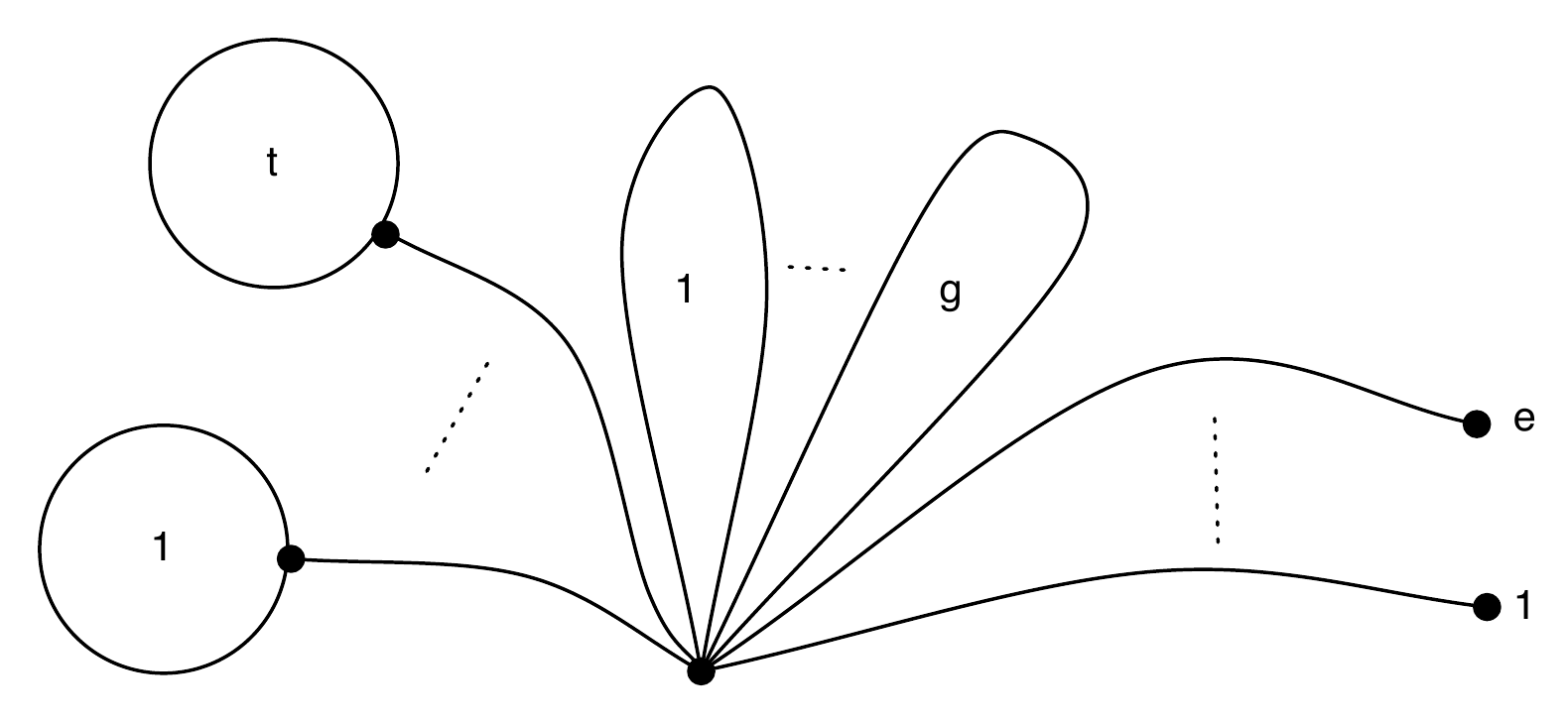}
\end{center}

The \emph{base vertex} $x$ of $G_v$ is the $0$-cell onto which the first edge
is attached. Let $\Hpty(G_v,\partial)$ be the space self-homotopy
equivalences of $G_v$ which,

\begin{enumerate}
\item fix the $e$ edges pointwise,
\item fix the $t$ loops of the boundary tori pointwise and
\item do not identify the base vertices of any two boundary tori.
\end{enumerate}

\begin{definition}{($H_v$)}\label{htyv}
Let $H_{v} = \pi_0\Hpty(G_{v},\partial)$ be the group of path components of
the space of self-homotopy equivalences described above. 
\end{definition}

When we write $v$ as $(g,i+o,a+b)$ we mean that the number of incoming edges
$i = |\In(G)|$, outgoing edges $o = |\Out(G)|$, incoming tori $a =
|\Tin(G)|$ and outgoing tori $b= |\Tout(G)|$.  If $[n]$ is the set
$\{1,\ldots,n\}$ then a \emph{boundary labelling} is a choice of
homeomorphisms, $i_H : [|\In(G)|]\times [0,1] \to \In(G)$ and $o_H :
[|\Out(G)|] \times [0,1] \to \Out(G)$. So that the interval ${i}\times
[0,1]$ is mapped homeomorphically onto the $i$th incoming or outgoing edge
and $i\times 0$ sent to the boundary vertex. For the tori we use the maps,
$$a_H : [|\Tin(G)|] \times [0,2\pi) \to \Tin(G)\quad\textnormal{ and }\quad b_H : [|\Tout(G)|] \times [0,2\pi) \to \Tout(G), $$

and we require that the points $a_H(i,0)$ and $b_H(i,0)$ are the base vertices
of the boundary torus. Compare to section \ref{gorprelim}.

\begin{definition}{(\cal{OCH})}\label{ochoh}
There is a symmetric monoidal category $\cal{OCH}$ enriched over $\Group$
with objects generated by the elements $e$ and $t$. The object $e^{\otimes
  n}$ represents $n$ labelled edges and the object $t^{\otimes k}$ represents $k$
boundary tori. The morphisms of $\cal{OCH}$ are self-homotopy equivalences
of boundary labelled graphs fixing boundary elements:

$$\Morph_{\mathcal{OCH}}(e^{\otimes i}\otimes t^{\otimes j}, e^{\otimes k}\otimes t^{\otimes l}) = \coprod_{g} H_{(g,i+k,j+l)}.$$

There are no morphisms between empty objects and we require $j\geq 1$ when
$l\geq 0$. The composition of $[\varphi] \in H_{v}$ and $[\psi] \in H_{w}$
is given by choosing maps $\varphi : G_v \to G_v$ and $\psi : G_w\to G_w$
which preserve the boundary labelling in $[\varphi]$ and $[\psi]$
respectively. 

\begin{center}\includegraphics{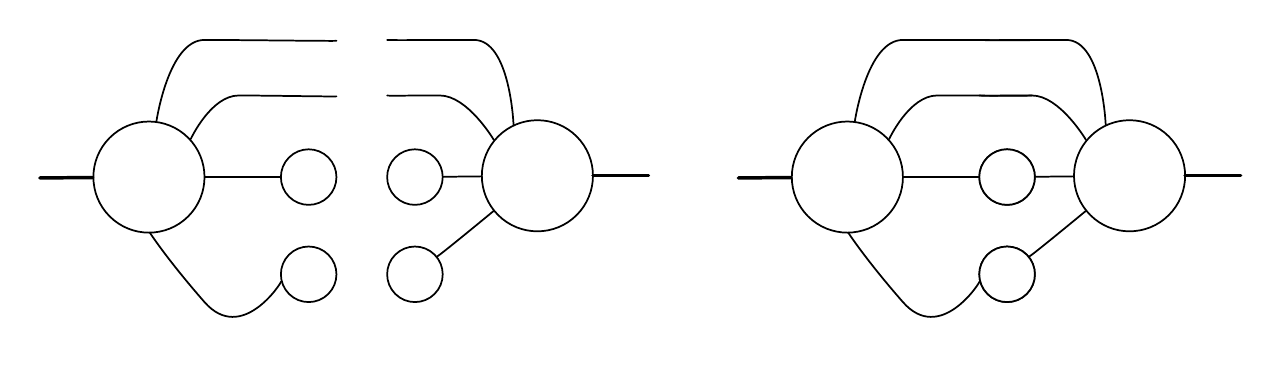}\end{center}

The illustration above depicts three graphs: $G_v$, $G_w$ and
$G_{v\#w}$. The graph $G_{v\#w}$ is formed by gluing together the graphs
$G_v$ and $G_w$. The dark lines represent boundary edges which are not
involved in the gluing. Using the standard embeddings found in section
\ref{hvwtheoremsec}, the picture above can be seen to correspond to the
gluing of manifolds $M_v$. 

Two homotopy equivalences can be glued to give an equivalence $\varphi\#\psi
: G_{v\#w} \to G_{v\#w}$. For any continuous variation of $\varphi$ or
$\psi$ within their respective path components, the graph $\varphi\#\psi$
varies continuously within the corresponding path component of
$\Hpty(G_{v\#w},\partial)$. The map $(\varphi,\psi)\mapsto \varphi\#\psi$
yields a composition law,
$$ H_{v} \times H_{w} \to H_{v\#w}.$$

This determines the composition law for the category $\cal{OCH}$.
\end{definition}

\begin{definition}{(\cal{OH})}\label{ohdef}
The \emph{open homotopy category}, $\mathcal{OH}$, is the subcategory of
$\mathcal{OCH}$ associated to graphs without tori.
\end{definition}

It follows from the discussion in \ref{dg categories} that there is a
monoidal category $B\mathcal{OCH}$ enriched over $\Top$. This category has
the same objects and its morphism spaces are equal to classifying spaces of
the groups defined above. Applying the functor $C_*(-;\mathbb{Q})$ yields
a differential graded category. 
%These categories will yield combinatorial
%analogues of the categories $\cal{O}$ and $\cal{OC}$ defined in section
%\ref{ocandoc}.

\begin{definition}{(\cal{OG}, \cal{OCG})}\label{opengr}
 The \emph{open graph category} $\cal{OG}$ and the \emph{open-closed graph
   category} $\cal{OCG}$ are the categories of rational chains on the
 classifying categories of the open and open-closed homotopy categories.
$$\cal{OG} =  C_*(B\cal{OH};\mathbb{Q}) \quad\normaltext{ and }\quad\cal{OCG} = C_*(B\cal{OCH};\mathbb{Q})$$
\end{definition}

\subsection{A Theorem of Hatcher, Vogtmann and Wahl} \label{hvwtheoremsec}

The theorem below appears in the papers of Hatcher, Vogtmann and Wahl.
It stems from Hatcher's work on the homotopy type of the diffeomorphism
group of $S^1 \times S^2$ \cite{MR1314940} and Vogtmann's study of Outer
Space \cite{MR830040}. The synthesis of these ideas has recently led to
homological stability results for 3-manifolds \cite{MR2113904,MR2174267}.

The mapping class groups in our construction will differ from those
considered in the references above by requiring that group elements fix a
regular neighborhood of the boundary (see section \ref{3dcob}). As such they
will be subgroups $\Gamma(M_v,\partial) \subset \Gamma(M_{v})$ generated by
the same generators given by Wahl and Jensen (\cite{MR2077676}) minus those
which require Dehn twists of the boundary torus. Differences will be noted
along the way.

\begin{definition}{($\Gamma_v$)}
The group $\Gamma_{v} = \Gamma(M_{v},\partial)$ is the mapping class group of the space $M_{v}$ considered in section \ref{ocandoc}. 
\end{definition}

Since $\pi_1(SO(3)) \cong \mathbb{Z}/2$, the inclusion $SO(3)
\hookrightarrow \Diff(S^2)$ yields a 1-parameter family of diffeomorphisms
$\varphi : S^2 \times I \to S^2$ such that one composition along the second
parameter is homotopic to identity. A \emph{Dehn twist} along a 2-sphere in
a 3-manifold is obtained by deleting a regular neighborhood of the sphere
and gluing the two boundary components back together along a copy of $S^2
\times I$ using the map $\varphi$.

We fix a \emph{standard embedding}, $i : G_v \hookrightarrow M_v$, by
mapping the end of each boundary edge $e$ to a boundary sphere, we require
each boundary torus of the graph to map to the loop on the longitude of
the boundary torus of $M_v$ and each of the $g$ loops to be sent to the $S^1$
component of the corresponding $S^1 \times S^2$ term. The inclusion $i$
induces an isomorphism on fundamental groups. Let $r : M_{v} \to G_{v}$ be
the retraction onto $i(G_{v})$ in $M_v$. These maps are canonical up to isotopy,
with respect to the decomposition of $M_{v}$ into punctured handle bodies.

{\Small
\begin{center}
\includegraphics[scale=.7]{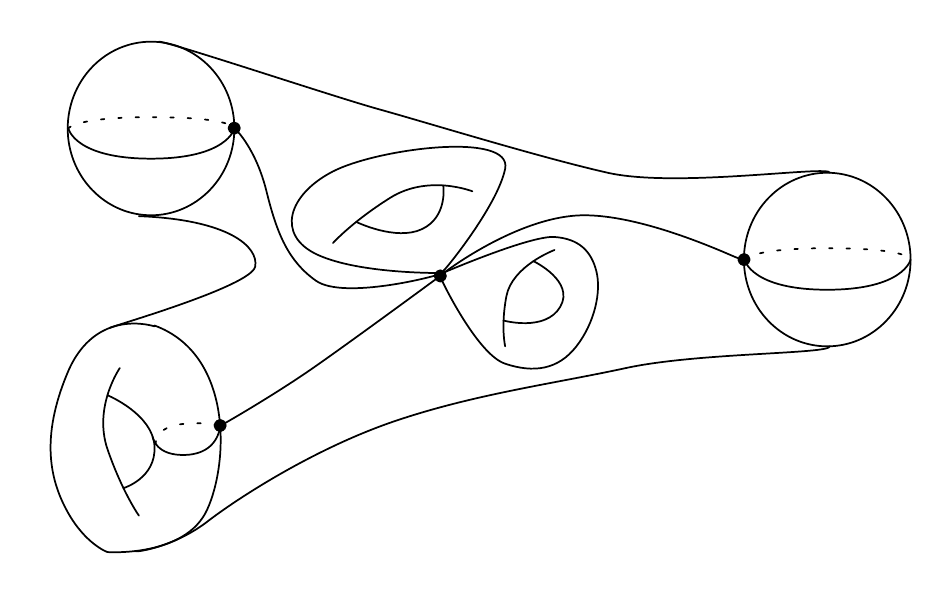}
$G_{(2,2,1)} \hookrightarrow M_{(2,2,1)}$
\end{center}
}

The illustration above consists of a graph $G_{(2,2,1)}$ embedded inside of
a 3-manifold $M_{(2,2,1)}$. The two internal loops inside of this graph
travel around the $S^1 \times S^2$ summands in the center of the
picture. The two boundary edges of the graph are attached to the two
boundary 2-spheres and the one boundary torus loop of the graph is attached
around the torus component of the 3-manifold.

If $l \in \Diff(M_{v},\partial)$ is a diffeomorphism, then we obtain a homotopy equivalence,
$$h(l) = r\circ l\circ i.$$

This defines a map $h : \Gamma_{v} \to H_{v}$. The key point for us is that
$h$ is a rational isomorphism, see corollary \ref{hratiso}. 

\begin{theorem}{(Hatcher-Vogtmann-Wahl)}\label{HVW}
  The map $h : \Gamma_{v} \to H_{v}$ is an epimorphism and its kernel is
  isomorphic to a finite direct sum of $\mathbb{Z}/2$'s generated by Dehn
  twists along spheres.
$$\begin{diagram}
 1 & \rTo & \bigoplus_k \mathbb{Z}/2 & \rTo & \Gamma_{v} &\rTo^h & H_{v} & \rTo 1
\end{diagram}$$
\end{theorem}

\begin{proof}
The reader may compare what follows to theorem 1.1 in \cite{MR2174267}.  In
their work Hatcher, Vogtmann and Wahl allow the mapping class groups above
to move the boundary while we do not. In our discussion of the difference,
we will simplify matters slightly by only discussing the tori.  If the
number of edges is equal to zero ($e=0$) then the full group of graph
automorphisms is generated by:

$$\begin{array}{rllrll}
1. & P_{i,j} & \textnormal{exchanges $x_i$ and $x_j$} & 5. & (x_i^{-1};y_j) & x_i \to y_j^{-1} x_i \\
2. & I_i & \textnormal{exchanges $x_i$ and $x_i^{-1}$} & 6. & (y_i^{\pm};x_j) & y_i \to x_j^{-1} y_i x_j \\

3. & (x_i;x_j) & x_i \to x_i x_j & 7. & (y_i^{\pm};y_j) & y_i \to y_j^{-1} y_i y_j. \\
4. & (x_i;y_j) & x_i \to x_i y_j & & & \\
\end{array}$$
%% $$\begin{array}{cll}
%% 1. & P_{i,j} & \textnormal{exchanges $x_i$ and $x_j$} \\
%% 2. & I_i & \textnormal{exchanges $x_i$ and $x_i^{-1}$} \\

%% 3. & (x_i;x_j) & x_i \to x_i x_j \\
%% 4. & (x_i;y_j) & x_i \to x_i y_j \\
%% 5. & (x_i^{-1};y_j) & x_i \to y_j^{-1} x_i \\

%% 6. & (y_i^{\pm};x_j) & y_i \to x_j^{-1} y_i x_j \\
%% 7. & (y_i^{\pm};y_j) & y_i \to y_j^{-1} y_i y_j. \\
%% \end{array}$$

The $x_i$ represent generators of $\pi_1(G_{(g,0,t)})$ associated to
factors of $S^1\times S^2$ and $y_i$ represent generators of
$\pi_1(G_{(g,0,t)})$ associated to factors of $S^1\times D^2$.

If we view our 3-manifold as the boundary of a punctured handle body then
generators 3-7 above can be represented by handle slides along the curves
$x_i$ and $y_j$. Handle slides are associated to generators of the
automorphism group as follows.

%% $$\begin{array}{cc}
%% 3. & \textnormal{The handle $x_i$ slides over $x_j$.}\\
%% 4. & \textnormal{The handle $x_i$ slides over $y_j$.}\\
%% 5. & \textnormal{The handle $x_i^{-1}$ slides over $y_j$.}\\

%% 6. & \textnormal{The torus $y_i$ slides over the handle $x_j$.}\\
%% 7. & \textnormal{The torus $y_i$ slides over the torus $y_j$.}\\
%% \end{array}$$
$$\begin{array}{llll}
3. & \textnormal{The handle $x_i$ slides over $x_j$.} & 6. & \textnormal{The torus $y_i$ slides over the handle $x_j$.}\\
4. & \textnormal{The handle $x_i$ slides over $y_j$.} & 7. & \textnormal{The torus $y_i$ slides over the torus $y_j$.}\\
5. & \textnormal{The handle $x_i^{-1}$ slides over $y_j$.} & & \\
\end{array}$$

In order to slide a handle or a torus (thought of as a connected sum of $S^1
\times D^2$'s) over a torus, a Dehn twist must be performed. Fixing the
boundary kills generators 4, 5 and 7. Since our homotopy groups are defined
to fix the loop of the graph contained in the torus, the correspondence is
preserved.
\end{proof}

\begin{corollary}\label{hratiso}
  The chain complexes $C_*(B\Gamma_{v};\mathbb{Q})$ and $C_*(BH_{v};\mathbb{Q})$ are quasi-isomorphic.
\end{corollary}

\begin{proof}
The map $Bh$ induces an equivalence because $B(\mathbb{Z}/2) \simeq
\mathbb{R}P^\infty$ and $\mathbb{R}P^{\infty}$ is rationally contractible.

  %% This follows because $B(\mathbb{Z}/2) \simeq \mathbb{R}P^\infty$ and
  %% $\mathbb{R}P^{\infty}$ is rationally contractible; so that $Bh$ induces an
  %% equivalence.
\end{proof}

The corollary above implies that the space of morphisms in the categories
$\cal{O}$ and $\cal{OC}$ (section \ref{ocandoc}) are rationally
quasi-isomorphic to those of $\cal{OG}$ and $\cal{OCG}$ respectively
(section \ref{opengr}). The theorem below follows from the observation that
the map inducing this equivalence is compatible with the gluing of open
boundaries.

\begin{theorem}\label{opencateq}
  The open category $\mathcal{O}$ of section \ref{ocandoc} is
  quasi-isomorphic to the open graph category $\mathcal{OG}$, see  definition \ref{opengr}.
$$\mathcal{O} \cong \cal{OG}$$
\end{theorem}
\begin{proof}
 The map $h$ as defined above is compatible with gluing the
 spherical boundary components,

$$\begin{diagram}
\Gamma_{v}\times \Gamma_{w} & \rTo^\# & \Gamma_{v\#w}\\
\dTo^h        &                   & \dTo^h\\
H_{v}\times H_{w} & \rTo^\# & H_{v\#w}\\
\end{diagram}$$

%where $v = (g,i+k,0)$, $w = (g',k+r,0)$ and $v\#w=(g+g'+k-1,i+r,0)$. Given
(see notation section \ref{subscriptnote}). Given $\varphi \in \Diff(M_v,\partial)$
 and $\psi \in \Diff(M_{w},\partial)$, the action of $\varphi\#\psi$ on
 $i(G_{v})\#i(G_{w})\subset M_{v\#w} = M_{v}\#M_{w}$ is the same as the
 action of $\varphi$ on $i(G_{v})$ glued to the incoming edges of $\psi$
 acting on $i(G_{w})$. This is because $\varphi$ and $\psi$ are required to fix a
 regular neighborhood of the boundary.

The maps $h$ induce a functor $\cal{O} \to \cal{OG}$. One can choose
sections of $h$, $H_{v} \to \Gamma_{v}$, so that there is a
functor $i : \cal{OG} \to \cal{O}$. We have $h\circ i = 1$ and $i\circ h
\simeq_{\mathbb{Q}} 1$.
\end{proof}

Recall the notion of the category $\Obj(\cal{D})$ associated to a
monoidal category $\cal{D}$ (see definition \ref{objcat} section
\ref{monoidal categories}).  The category $\cal{OC}$ defines an
$\Obj(\cal{OC})-\cal{O}$ bimodule,
$$\cal{OC} : \Obj(\cal{OC})\otimes\cal{O}^{op}\to \Chain$$

via $(e^{\otimes n}\otimes t^{\otimes m})\otimes o^{\otimes k} \mapsto
\Morph(o^{\otimes k}, e^{\otimes n}\otimes t^{\otimes m})$, see also the
observation in section \ref{dg modules}.

The category $\mathcal{OCG}$ (definition \ref{opengr}) defines an
$\Obj(\cal{OCG})-\cal{OG}^{op}$ bimodule in the same way.  In fact, the
category $\cal{OCG}$ also defines an $\Obj(\cal{OC})-\cal{O}$ bimodule
because $\Obj(\cal{OCG}) = \Obj(\cal{OC})$ and theorem \ref{opencateq} above
implies that $\cal{O} \cong \cal{OG}$.

We have two $\Obj(\cal{OC})-\cal{O}$ bimodules: $\cal{OC}$ and $\cal{OG}$.
Corollary \ref{hratiso} shows that these two bimodules are
the same.

\begin{theorem}\label{octoocg}
  As $\Obj(\cal{OC})-\cal{O}^{op}$ bimodules the categories $\cal{OC}$ and
  $\cal{OCG}$ are quasi-isomorphic.
\end{theorem}

\subsection{Outer Space} \label{opouterspacesec}
For each $v = (g,e,t)$, we begin by defining a set $L_v$ consisting of
labelled graphs. This set will be used to construct a simplicial set
\Lp. The geometric realization \Yp\ of \Lp\ will be a classifying space for the
group $H_v$. In what follows, all graphs will be boundary labelled and we
will consistently write $L_v$ where $v = (g, i+o, a+b)$, see section \ref{hty
  eq grp}.

A graph $G$ is \emph{labelled} when it is paired with a map $\varphi :
G_{v}\to G$ which satisfies the following properties.

\begin{enumerate}
\item The function $\varphi$ preserves the incoming and outgoing edges and
  identifies the ends of each of the boundary tori of $G_{v}$ with circles
  $G$. By circle we mean cycles with one edge and one vertex.

\item If $x$ is the vertex of $G_{v}$ then the induced map, $\varphi_* :
  \pi_1(G_{v},x) \to \pi_1(G,\varphi(x))$ is an isomorphism.
\end{enumerate}

Two labelled graphs $(G,\varphi)$ and $(G',\psi)$ are \emph{equivalent} if
there is a graph isomorphism $\rho : G \to G'$ so that the diagram below
commutes.

$$\begin{diagram}
\pi_1(G,\varphi(x)) &             & \rTo^{\rho_*}              &       & \pi_1(G',\psi(x))\\
            & \luTo^{\varphi_*}        &                   & \ruTo^{\psi_*} & \\
            &              & \pi_1(G_{v},x) &    &\\
\end{diagram}$$

\begin{definition}{($L_v$)} \label{opensimp}
If $v= (g,e,t)$ then $L_v$ will denote the set of equivalence classes
$(G,G_v\xrightarrow{\varphi} G)$ of labelled graphs.
\end{definition}

The set $L_{v}$ can be endowed with a simplicial structure in which the
faces of simplices are determined by edge collapses (see section
\ref{gorprelim}). In what follows, we will use the nerve \Lp\ of $L_v$. A non-degenerate $n$-simplex in \Lp\ is given by a sequence
$$(G_0,\varphi_0) \subset (G_1,\varphi_1) \subset \cdots \subset (G_n,\varphi_n)$$

where $(G_i,\varphi_i) \in L_v$ for each $i=0,1,\ldots,n$ and $(G_i,\varphi_i)$
is obtained from $(G_{i+1},\varphi_{i+1})$ by collapsing one or more edges
(while preserving the homotopy type). Equivalently, simplices of the space
\Lp\ are determined by fixing a forest, $F_0 \subset G$, and a nested
sequence of subforests, $ F_n \subset F_{n-1} \subset \cdots \subset F_0
\subset G$. If $\varphi$ is a labelling of $G = G_n$ then this gives the
simplex,
$$(G/F_0,\bar{\varphi}_0) \subset (G/F_{1},\bar{\varphi}_{1}) \subset \cdots \subset (G_n/F_{n},\bar{\varphi}_n).$$ 

The maps, $\bar{\varphi}_i$, are induced by collapsing edges. In what follows we
will require all forests $F \subset G$ to
\begin{enumerate}
\item include all of the vertices of $G$,
\item include \emph{none} of the incoming or outgoing open boundary edges and
\item include no two base vertices of tori in the same tree.
\end{enumerate}

Simplicial face maps are defined by combining collapses and simplicial
degeneracy maps are given by inserting identity collapses.

The group $H_v$ acts on the nerve \Lp\ by changing the labellings. If $f \in H_{v}$ then $f : L_{v}
\to L_{v}$ is defined by $f(G,\varphi) = (G,\varphi\circ f)$ and so $f :
\Lp \to \Lp$ acts by
\begin{align*}
(G/F_0,\varphi_0) \subset (G/F_1,\varphi_1) \subset & \cdots \subset (G/F_n,\varphi_n)\\  \mapsto & (G/F_0,\varphi_0\circ f) \subset (G/F_1,\varphi_1\circ f) \subset \cdots  \subset (G/F_n,\varphi_n\circ f).
\end{align*}

\begin{definition}{(\Lp, \Yp, \Xp)}
The geometric realization of \Lp\ will be denoted by \Yp\ and the quotient
$\Yp/H_{v}$ will be denoted by either \Xp\ or \Xpq, see theorem \ref{modthm}
below.
\end{definition}

Suppose that $v = (g,e,t)$, if $t=0$ and $e = 0$ then \Xp\ is called Outer
space since the construction is a model for the classifying space of the
group of outer automorphisms of the free group $F_g$, see
\cite{MR830040}. If $t = 0$ and $e = 1$ then \Xp\ is known as ``Auter
space.''  Other generalizations, not involving diffeomorphisms that fix the
boundary, can be found in \cite{MR2113904,MR2077676, MR2174267}.

\begin{theorem}{($\Yp/H_v$ models $BH_v$)}\label{modthm}
  The action of $H_{v} = \pi_0\Hpty(G_{v},\partial)$ on the space
  \Yp\ is properly discontinuous and the stabilizer of any given
  simplex is a finite group. Moreover, the space \Yp\ is contractible.
\end{theorem}

\begin{proof}
 The action of $H_v$ is almost free. If $f\in H_{v}$ then $f(G,\varphi) =
 (G,f\circ \varphi) = (G,\varphi)$ if and only if $f$ is an isomorphism
 of the graph $G$. A graph isomorphism is determined by the manner in
 which it permutes the edges and so the size of the group of graph
 isomorphisms is bounded above by the group of all permutations on edges.

The proof of contractibility of \Yp\ is a special case of the proof which
appears in Wahl and Jensen's article \cite{MR2077676}.
\end{proof}

\begin{corollary}
  The quotient space $\Xp = \Yp/H_{v}$ is a rational model for the
  classifying space of $H_{v}$. In particular,
$$C_*(BH_{v};\mathbb{Q}) \simeq C_*(\Xp; \mathbb{Q}). $$
\end{corollary}

There is a geometric interpretation of the space \Xp. A \emph{metric
  graph} is a graph together with a fixed length $l(e) \geq 0$ assigned to
each internal edge. A metric graph is \emph{balanced} if $\sum_{e\in E(G)}
l(e) = 1$. The space \Xp\ is a subdivision of the space of balanced metric
graphs homotopy equivalent to the graph $G_{v}$.  For any balanced metric
graph $G$, if $e_0,\ldots, e_k$ are its edges then $G$ is uniquely
represented by the barycentric coordinates
$(l(e_0),\ldots,l(e_k))$ of a $k$ simplex $\Delta$ associated to the
topological type of $G$.

The boundary tori are represented by balloons attached to the graphs
representing points in the moduli space \Xp. The length of the edge at
the end of each balloon is fixed. The length of the edge used to attach the
balloon to the rest of the graph is allowed to vary and may
approach zero providing two distinct base vertices do not touch as a result.

We metrize the graphs in this way because the edge of the balloon
corresponding to a torus in a manifold $M_{v}$ is completely fixed by the
action of any $b\in \Gamma(M_{v},\partial)$. The edge about the torus
in the graph $G_{v}$, thought of as embedded in $M_{v}$, does not vary
with respect to the action of the mapping class group. The edge that is used
to attach the balloon to the rest of the graph \emph{is} allowed to vary
because $b$ may move the boundary torus about inside of $M_{v}$. Since there
are disjoint regular neighborhoods of the boundary tori in the construction
of the cobordism category, we can ask for the base vertices of the balloons
representing them not to touch.

In contrast, the open edges are given fixed length. When represented as a
graph within $M_{v}$, this length reflects the disjointness of the regular
neighborhoods of $2$-spheres in the construction of the cobordism
category. Allowing these lengths to vary is not necessary and would not add
anything to what follows. If we allowed the lengths to vary then it would be
necessary for us to consider the scenario in which the collapse of an edge
represented a boundary collision as we have done with the tori above.

%or equivalently not given
%length in the moduli space. If taken to have a fixed positive length then

\subsection{Cellular Stratification by Cubes}\label{cubestrat}

In order to compute the homology of $X_v$, we group simplices that can be
obtained from the same forest into a single cell (see \cite{MR1671188,
  MR1341841, MR2026331}). The cells obtained from this construction will be
called cubes.

A cube $[G,F,\varphi] \subset \Yp$ is obtained by gluing together all the
simplices arising from different filtrations of some fixed forest
$F\subset G$.
%% for a given labelled graph $(G, \varphi) \in L_{v}$.

$$ [G,F,\varphi] = \coprod_{F_0 \subset \cdots \subset F_m \subset F} (G/F_0 \subset \cdots \subset G/F_{m-1} \subset G/F_m) \times \Delta^m $$

The collection of all cubes $[G,F,\varphi]$ gives \Yp\ the structure of a CW
complex called the \emph{forested graph} stratification. Each cube
$[G,F,\varphi]$ in \Yp\ is homeomorphic to a $k$-cube $[0,1]^k$, where $k =
|E(F)|$. One can define such a homeomorphism by assigning an axis to each
edge.

If the graphs $G$ are planar trees then a construction analogous to the one
in section \ref{opouterspacesec} produces simplicial subdivisions of
associahedra. The cubical stratification above yields the cubical
decomposition of associahedra in this context, see \cite{MR0420609}.

The codimension 1 faces of a cube $[G,F,\varphi]$ are given by two
operations on graphs.

\begin{enumerate}
\item Collapsing an edge. $[G,F,\varphi] \mapsto
  [G/e,F/e,\bar{\varphi}]$ for some edge $e \in E(F)$.
\item Removing an edge from the forest. $[G,F,\varphi] \mapsto
  [G,F-e,\varphi]$ for some edge $e \in E(F)$.
\end{enumerate}

\begin{figure}\label{difffig}
\begin{diagram}
\BPic{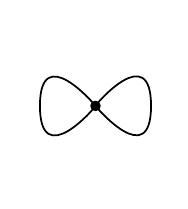} & \lTo^{(1)} & \BPic{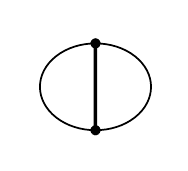} & \rTo^{(2)} & \BPic{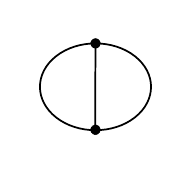} \\
\end{diagram}
\caption{The two faces of a cube $(G,F)$. The theta graph $G$ is drawn using light edges. Dark edges are used represent the forest $F\subset G$.}
\end{figure}

The two types of faces, $(1)$ and $(2)$, are illustrated in the figure. The
group $H_{v}$ now acts cellularly. The stabilizer of the cube
$[G,F,\varphi]$ consists of automorphisms of $G$ that send the forest $F
\subset G$ to itself.

Each cube $[G,F,\varphi]$ in \Yp\ descends to a cube $[G,F]$ in the
quotient \Xp. This cube is not necessarily a cell, but an
orbi-cell. This follows from identifying the cube in \Yp\ with a cube $C
= [0,1]^k$ where each edge of $F$ is associated to an axis. The portion of
the cube that descends to \Xp\ is the quotient of $C$ by the stabilizer
$\Aut(G,F,\varphi)$.  The action of $\Aut(G,F,\varphi)$ on $C$ fixes the
origin and permutes the axes so that $C/\Aut(G,F,\varphi)$ is a cone on the
quotient of the boundary $\partial C$.

\begin{lemma}\label{quotientlemma}
  The quotient of an $n$-sphere by a finite linear group $G\subset
  GL_n(\mathbb{R})$ is $\mathbb{Q}$-homotopic to either a $n$-sphere or a
  $n$-ball. The latter case holds only when the action includes reflections.
\end{lemma}

For proof and discussion, see \cite{MR1671188}. Those cubes which have
symmetries that do not include reflections survive to the quotient.

In \Xp\ the tori are represented by trees containing the base vertex of the
balloons.

\subsection{Homology}\label{homology}

In this section we complete our description of the morphism spaces of
$\cal{OC}$ and $\cal{O}$.  For each $v=(g,e,t)$, we define a generalized
$\Cobar$ construction: an exact functor $\cal{G}_{v}$ from the category of
differential graded cooperads to chain complexes. The complexes $\cal{G}_v$
will be those that generate the morphism spaces of the enveloping functor
$\Cobar(\cal{O})^{\flat}$ defined in \ref{rel to dg alg}. We will show that
$\cal{G}_v$ corresponds to the chain complex obtained from the
stratification of \Xp\ by cubes defined in the previous section.

\subsubsection{From Operads to Graph Complexes}\label{optogcplx}

A \emph{bonnet} is a graph $B(n)$ isomorphic to a corolla with two edges
identified.

$$B(n) =  \BPic{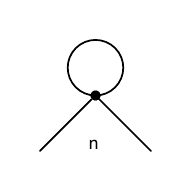}$$

Let $S_{v}$ be the set of boundary labelled combinatorial graphs of genus
$g+t$ with $e$ boundary edges and $t$ bonnets. A graph $G\in S_v$ differs
from $G_v$, pictured in section \ref{hty eq grp}, in that $G$ is allowed to
possess internal edges of any kind. Given a cyclic dg cooperad $\cal{P}$,
\emph{the generalized $\Cobar$ construction}, $\cal{G}_v(\cal{P})$, is the
complex consisting of $S_v$ graphs labelled by $\cal{P}$ and oriented using
the convention described in section \ref{orprelim}.

\begin{definition}{($\cal{G}_v(\cal{P})$)}

$$\cal{G}_{v}(\cal{P}) = \bigoplus_{G\in S_v} \cal{P}(G)\otimes\det(G)^*$$

\end{definition}

The differential $\delta$ expands edges. It can be described by its matrix
elements, $(\delta)_{G',G}$, where $G,G'\in S_v$.  If $G'$ is not isomorphic
to $G/e$ for some collapsible edge $e\in G$ then set
$(\delta)_{G',G}=0$. Otherwise, let $c : G \to G' \cong G/e$ so that if $c^*
: \mathcal{P}(G') \to \mathcal{P}(G)$ is the induced map on the labelling
then $\delta$ is given by $ (\delta)_{G',G} = c^* \otimes p^e$ where $p^e$
is the map induced on the orientation by collapsing the edge. If the cooperad
$\cal{P}$ has a non-trivial differential then the total differential is the
sum of the differential defined above together with the original
differential.

The generalized $\Cobar$ construction is introduced in order to mediate
between the algebraic world of operads and categories, and the topological
world obtained the from moduli spaces defined in earlier sections. In
particular, the collection $\{\cal{G}_v\}$ naturally models the morphisms of
the open and open-closed categories introduced in section \ref{ocandoc}. By
construction, we have the following identifications,

$$\cal{G}_{(0,e,0)}(\cal{P}) = \Cobar(\cal{P})(e)\quad\normaltext{and}\quad\Morph_{\Cobar(\cal{P})^\flat}(x^{\otimes n},x^{\otimes m}) = \bigoplus_g \cal{G}_{(g,n+m,0)}(\cal{P}).$$

\begin{remark}\label{remyadda}
We can use this observation to relate the $\cal{G}_v$ to modular operads.
In particular, when $t=0$ the collection $\{ \cal{G}_v \}$ determine a PROP,
$\Cobar(\cal{P})^\flat$, see section \ref{rel to dg alg}. By remark
\ref{reltomodularremark}, the PROP $\Cobar(\cal{P})^\flat$ agrees with
$P\cal{M}\Cobar(\cal{P})$.  On the other hand, the Feynman transform (see
\cite{GetzKapMod}) of the modular operad associated to a cyclic operad
commutes with the $\Cobar$ construction (with appropriate twisting),
$\cal{M}\Cobar(\cal{O}^\vee) \cong {\mathsf F} \cal{M}\cal{O}$. This yields a
relationship between the generalized $\Cobar$ construction and the Feynman
transform, $\bigoplus_{t=0}\cal{G}_v(\cal{P}) \cong
P{\mathsf F}\cal{M}(\cal{P}^\vee)$.
\end{remark}

If $t > 0$ then it is best to think of the collection $\{\cal{G}_v\}$ as
describing the extension of the $t=0$ case by data coming from the torus
boundary; a dg module over the open category. This comment is made more
precise in section \ref{extandtorussec}.

\begin{lemma}\label{exactnesslemma}
The functor $\cal{G}_v$ is exact: if $\varphi : \cal{P} \to \cal{P}'$ is a
  quasi-isomorphism of cooperads then the induced map $\cal{G}_v(\cal{P})\to
  \cal{G}_v(\cal{P}')$ is a quasi-isomorphism.
\end{lemma}

This is proven using a spectral sequence argument, see
\cite{GetzKapMod} Theorem 5.2 (3).

\subsubsection{Cubical Chains Compute A Double Dual}

Recall from section \ref{cubestrat} that the complex of cubical chains on
\Xp\ is spanned by cubes $[G,F]$ where $G$ is a boundary labelled graph
with $t$ cycles representing boundary tori and $F \subset G$ is a forest
containing all of the vertices of $G$ and none of the boundary edges. No two
vertices of the boundary tori are contained in the same tree of $F$.

The cube $[G,F]$ is oriented by an ordering of the edges of $F$. Lemma
\ref{quotientlemma} in the same section implies that the antisymmetry
relation $[G,-F] = -[G,F]$ holds.

The differential is given by the sum over ways to remove an edge from a
forest and the sum over ways to contract an edge contained in the forest. In
either case the cube is oriented by the induced orientation.

$$ \partial [G,F] = \sum_{e\in F} [G/e,F/e] + \sum_{e\in F} [G,F-e]$$

Recall that $\Comm$ is the commutative operad defined in section
\ref{examples of operads}.  The cooperad $\BBar(\Comm)$ is the free cooperad
on $n$-corolla satisfying the antisymmetry relation (dual to the
$L_{\infty}$ operad). The trees are edge oriented and the differential
contracts edges.

Since $\Cobar(\BBar(\Comm))$ is a double complex, while $C_*(\Xp)$ is merely a
chain complex, we flatten the double grading as follows,

$$ \Cobar(\BBar(\Comm))(n)^{'}_i = \bigoplus_j \Cobar(\BBar(\Comm))(n)_{j,i} .$$

The differential $d$ remains the sum of the internal differential, which
contracts the edges of $\BBar(\Comm)$, and the external differential, which
expands compositions.

\begin{theorem}\label{foresttocobar}
  The rational homology of the spaces \Xp\ is computed by
  $\cal{G}_v(\BBar(\Comm))$:
  $$C^{cell}_*(\Xp;\mathbb{Q}) \cong \cal{G}_{v}(\BBar(\Comm))' .$$
\end{theorem}

\begin{proof}
  Assuming $t=0$, by lemma \ref{envelopingequiv} it suffices to show that the
  operad $\Cobar(\BBar(\Comm))$ is isomorphic to the operad with $\cal{O}(n)
  = C_*(\Xr_{(0,n+1,0)};\mathbb{Q})$. This forms an operad because the
  cellular composition, theorem \ref{opengluing}, is independent of this
  theorem.  We will see that as complexes the two are plainly isomorphic:
  $$ \Cobar(\BBar(\Comm))(n)' \cong C_*(\Xr_{(0,n+1,0)};\mathbb{Q}).$$

  In degree $j$, the complex $C^{cell}_j(\Xr_{(0,n+1,0)};\mathbb{Q})$ is
  spanned by forested trees, $(T,F)$, where the forest $F$ contains $j$ edges
  and a connected component associated to each internal vertex of $T$. 

  In bidegree $(j,i)$, the complex $\Cobar(\BBar(\Comm))(n)_{j,i}$ is
  spanned by unrooted $n$ trees $T$, containing $j = |T|$ internal vertices
  each of which is labelled by a tree $F_l \in \BBar(\Comm)(H(v))$. The
  equation $(j,i) = (|T|, \sum_{m=1}^{|T|} (|F_m| - 1))$ holds for
  bidegrees. Since the second coordinate is the total number of internal
  edges, $T \otimes F_1 \otimes \cdots \otimes F_j \in
  \Cobar(\BBar(\Comm))(n)^{'}_i$ when $T$ is an unrooted $n$ tree labelled by
  trees $F_l$ whose internal edges total to $i$.

  %% if we include operad suspension
  %% In degree $n$ the complex $C^{cell}_j(X'_{(0,n+1,0)};\mathbb{Q})$ is spanned by
  %% forested trees $[T,F]$ where the forest $F$ contains $j$ edges and a
  %% connected component associated to each internal vertex of $T$.

  %% In bidegree $(i,j)$ the complex $\Cobar(\BBar(\Comm))(n)_{i,j}$ is spanned
  %% by trees $T$ containing $-j = |T|$ internal vertices each of which is in
  %% turn labelled by a tree $F_l \in \BBar(\Comm)(H(v))$ the degree of $F_l$
  %% is $|F_l|$. So that the total degree is $n = i+j = -|T|+\sum_{l=1}^{|T|}
  %% |F_l|$. If $E(F_l)$ is the number of internal edges of $F_l$ we have
  %% $|E(F_l)| + 1 = |F_l|$ which implies $n = \sum_{m=1}^{|T|} |E(F_l)|$.
  %% Thus for $T \otimes F_1 \otimes \cdots \otimes F_j \in
  %% \Tot\Cobar(\BBar(\Comm))(n)_n$ if $T$ is tree labelled by trees $F_l$
  %% whose internal edges total to $n$.
  
  To a forested tree $[T,F]$ with $F = F_1 \cup \cdots \cup F_j$, we
  associate the tree with internal vertices labelled by the $F_l$. The
  inverse map is obtained by doing the opposite: inserting forests at
  vertices.

  The two differentials in either complex are the same. Collapsing an edge
  in a forest corresponds to contracting an edge in a $\BBar(\Comm)$
  labelling. Removing an edge in a forest corresponds to inserting an edge
  in $\cal{G}_v$ between two $\BBar(\Comm)$ labellings; this is the $\Cobar$
  differential. See figure 1 in section \ref{cubestrat}.

  The two orientation conventions agree. A forested graph $[T,F]$ is
  oriented by an ordering of the edges in the forest $F$. If $F = \cup_i
  F_i$ then

$$\det(E(F)) = \bigotimes_i \det(E(F_i)).$$

On the other hand, if a graph $G$ is a tree $T$ with $j$ vertices labelled by
forest components $F_1,\ldots, F_j$ then, the convention described in section
\ref{orprelim} tells us that,

$$\det(T\otimes F_1\otimes \cdots\otimes F_j) = \det(E(T))\otimes \det(\Out(T)) \bigotimes_{i=1}^j \det(E(F_i))\otimes \det(\Out(F_i)).$$

In our computation, the number of outgoing edges of $T$ is one.  The
internal edges of $T$ join the labellings of two separate vertices by forest
components $F_i$. One end of each edge of $T$ is an incoming edge of some
forest component and the other end is an outgoing edge of some forest
component.

The outgoing components of each forest must correspond to internal edges of
$T$ except for the one outgoing edge corresponding to the outgoing edge of
$T$. Thus there is a bijection between the set $E(T)\coprod \Out(T)$ and
$\coprod_i \Out(F_i)$. Taking graded determinants yields the isomorphism,

$$\det(E(T)) \cong \det(E(T))\otimes \mathbb{Q} \cong \bigotimes_{i=1}^j \det(\Out(F_i)).$$

It follows that $\det(T\otimes F_1\otimes \cdots\otimes F_j) \cong \otimes_i
\det(E(F_i))$ and so the signs in both differentials agree.

If the number of tori is greater than zero then the cells associated to the
boundary tori are the trees containing the base vertex of the balloon
associated to the torus. These are represented combinatorially by bonnets in
$\cal{G}_v(\BBar(C))$.
\end{proof}

%\subsection{Corollaries}

\begin{corollary}\label{reductioncor}
$$\Morph_{\mathcal{O}}(e^{\otimes i}, e^{\otimes j}) \simeq \Morph_{\Cobar(\BBar(C))^\flat}(e^{\otimes i}, e^{\otimes j})$$
\end{corollary}

Recall the notion of the category $\Obj(\cal{D})$ associated to a monoidal
category $\cal{D}$ (definition \ref{objcat} section \ref{monoidal
  categories}).

\begin{corollary}\label{reductioncor2}
  The $\Obj(\cal{OC})-\cal{O}$ bimodule $\mathcal{OC}$ is quasi-isomorphic to
  the $\Obj(\cal{OC})-\cal{O}$ bimodule defined by the functor,

$$ (e^{\otimes n}\otimes t^{\otimes m})\otimes o^{\otimes k} \mapsto \coprod_g \cal{G}_{(g,n+k,m)}(\BBar(C)).$$
\end{corollary}

 The corollary follows from the identification, $ C_*(BH_{v};\mathbb{Q})\simeq C^{cell}_*(\Xp;\mathbb{Q}) $ and the previous theorem.

\section{The Open Category}\label{opencatgluing}

Corollary \ref{reductioncor} states that morphisms of the category $\cal{O}$
are quasi-isomorphic to spaces of graphs. In this section we show that the
composition induced from the gluing of 2-spheres in the open category is
cellular. This allows us to extend corollary \ref{reductioncor} from an
equivalence of morphism spaces to an equivalence of categories.  The
combinatorial open category $\Cobar(\BBar(C))^{\flat}$ is equivalent to the
open category $\cal{O}$.

Given two boundary labelled composable forested graphs $[G,F]$ and
$[G',F']$, form the graph $G\#G'$ by gluing the relevant ends together and
eliminating the resulting bivalent vertices. The forests $F$ and $F'$
together form a forest $F\cup F'$ of $G\#G'$, because forests are not
permitted to contain boundary edges.

\begin{theorem}\label{opengluing}
  The quasi-isomorphisms of \ref{reductioncor} respect composition.

{\Small $$\begin{diagram}
\Morph_{\mathcal{OG}}(e^{\otimes i}, e^{\otimes j}) \otimes \Morph_{\mathcal{OG}}(e^{\otimes j}, e^{\otimes k})
 & \rTo^{\circ} & \Morph_{\mathcal{OG}}(e^{\otimes i}, e^{\otimes k}) \\
\dTo^{\varphi_{ij}\otimes \varphi_{jk}}        &                   & \dTo^{\varphi_{ik}}\\
\Morph_{\Cobar(\BBar(\Comm))^{\flat}}(e^{\otimes i}, e^{\otimes j}) \otimes \Morph_{\Cobar(\BBar(\Comm))^{\flat}}(e^{\otimes j}, e^{\otimes k})
 & \rTo^{\circ} & \Morph_{\Cobar(\BBar(\Comm))^{\flat}}(e^{\otimes i}, e^{\otimes k}) \\
\end{diagram}$$ }

\end{theorem}

\begin{proof}
We show that composition respects the cube decomposition of the outer
spaces. In everything to follow, whenever the subscript $v = (g,e,t)$ is
used, we will assume that $t=0$.

The composition of $\mathcal{OG}$ is defined by maps,
$$\circ : C_*(X_v;\mathbb{Q}) \otimes C_*(X_{w};\mathbb{Q}) \to C_*(X_{v\#w};\mathbb{Q}).$$
%%$$\circ : C_*(BH_v;\mathbb{Q}) \otimes C_*(BH_{w};\mathbb{Q}) \to C_*(BH_{v\#w};\mathbb{Q}).$$

There are $\mathbb{Q}$-homotopy equivalences from the space $BH_{v}$ to $\Xp
= \Yp/H_v$. These spaces are stratified by orbi-cells, $[G,F]$, indexed by
forested graphs. The dimension of $[G,F]$ is the number of edges in $F$.
Residing above each orbi-cell is a collection of honest cells,
$[G,F,\varphi]$ in $\Yp$, which are indexed in the orbit of the action of
$H_{v}$ by their labellings $\varphi$, see \ref{cubestrat}.

Given a cell $[G,F,\varphi]$ of dimension $n$ in \Yp\ and a cell
$[G',F',\varphi']$ of dimension $m$ in \Yq\ (representing a pair of
composable graphs) there is a composite $[G\#G', F\cup
  F',\varphi\#\varphi']$ of dimension $n+m$ and a homeomorphism,
$$ [G,F,\varphi] \times [G',F',\varphi'] \to [G\#G', F\cup F',\varphi\#\varphi']$$

defined by identifying each cell with a cube in $\mathbb{R}^{|E(F)|}$ as
described in \ref{cubestrat}.  These homeomorphisms together yield a
composition,
$$ \Yp \times \Yq \to \Ypq $$

which is equivariant with respect to the action of $H_{v}\times H_{w}$ on
the left and $H_{v\#w}$ on the right, using the map $H_{v}\times H_{w} \to
H_{v\#w}$ (see definition \ref{ochoh}). So there is a composition on the
quotient. The composition of two cubes $[G,F]$ and $[G',F']$ is determined
by the diagram below.

\begin{center}$\begin{diagram}
[G,F,\varphi] \times [G',F',\varphi'] & \rTo^{\cong} & [G\#G',F\cup F',\varphi\#\varphi']\\
\dTo        &                   & \dTo\\
[G,F] \times [G',F'] & \rTo & [G\#G',F\cup F']\\
\end{diagram}$\end{center}

It can be seen that the differential acts as a derivation with respect to
this composition law using the rule in section \ref{cubestrat}.
\end{proof}

The theorem above, together with theorem \ref{foresttocobar}, implies the
following corollary.

\begin{corollary}\label{omodarecinfty}
  The category of h-split $\mathcal{O}$ modules is equivalent to the
  category of cyclic $\Cobar(\BBar(\Comm))$ algebras. In particular, the
  category of h-split $\mathcal{O}$ modules is equivalent to the category of
  cyclic $\Comm_{\infty}$ algebras.
\end{corollary}

It is possible to restate the result of theorems \ref{opengluing} and
\ref{foresttocobar} in the language of cyclic operads. Let $M_n = \#^n D^3$
be the 3-manifold obtained by connect summing $n$ copies of the 3-ball,
$D^3$, to itself. If we set
$$\mathcal{H}_n = C_*(B\Gamma(M_n,\partial);\mathbb{Q})$$ then the
collection $\{\mathcal{H}_n\}$ form a cyclic dg operad $\mathcal{H}$ quasi-isomorphic to
$\Cobar(\BBar(\Comm))$ where $\Comm$ is the commutative operad. The
machinery of modular operads implies the following corollary, see
\cite{GetzKapMod}.

\begin{corollary}\label{modularrefcor}
Cyclic $\Comm_{\infty}$ algebras are algebras over the modular closure of the chain operad $\mathcal{H}$ defined above.
\end{corollary}

\section{Extension and the Torus}\label{extandtorussec}

Given a cyclic $\Comm_{\infty}$ algebra $A$, corollary \ref{omodarecinfty}
shows that $A$ defines an open TFT in the sense of definition
\ref{octft}. From section \ref{dg modules}, the inclusion $i : \cal{O} \to
\cal{OC}$ induces a derived pushforward,
$$\mathbb{L}i_* : \cal{O}\module \to \cal{OC}\module.$$

Thus any such algebra $A$ determines an open-closed topological field theory
$\mathbb{L}i_*(A)$. On the other hand, the inclusion $j : \cal{C} \to
\cal{OC}$ determines a closed TFT, $j^* \mathbb{L}i_*(A)$. A closed TFT is a
$\cal{C}$ module. The $\cal{C}\module$ structure on $j^* \mathbb{L}i_*(A)$
is equivalent to the existence of a natural map,
$$\cal{C}(t^{\otimes i}, t^{\otimes j}) \otimes j^*\mathbb{L}i_*(A)(t^{\otimes i}) \to j^*\mathbb{L}i_*(A)(t^{\otimes j}).$$

In this section, we show that the homology of the chain complex associated
to the torus object, $j^*\mathbb{L}i_*(t)$, is the Harrison homology of the
algebra $A$. This is proven by studying the $\Obj(\cal{OC})-\cal{O}$
bimodule $\cal{OC}$ used to define the extension $\mathbb{L}i_*$ above.

Recall that the boundary tori in the forested graph stratification of the
space \Xp\ are represented by \emph{bonnets}, $B(n)$, see section
\ref{optogcplx}. The boundary of the trivial bonnet, $B(0)$, is zero.  In
general, the boundary of the cell associated to the tori derives from the
differential in the $\Cobar$ construction.

%% In what follows we see that the bimodule $\mathcal{OC}$ consists of graphs
%% from the category $\mathcal{O}$ decorated by bonnets. In particular, the
%% bonnets generate $\cal{OC}$ as a $\Obj(\mathcal{OC})-\mathcal{O}$ bimodule.

\begin{theorem}\label{bonnetgen}
  The category $\cal{OC}$, when considered as an
  $\Obj(\mathcal{OC})-\mathcal{O}$ bimodule, is freely generated by the
  bonnets $B(n)$.
\end{theorem}

\begin{proof}
It follows from corollary \ref{reductioncor2} that we can consider
$\cal{G}_v(\BBar(\Comm))$. If $G\in \cal{G}_{(g,n+k,m)}(\BBar(\Comm))$ is a
basis element then $G$ is a $\BBar(\Comm)$ labelled graph with $n$ incoming
edges, $k$ outgoing edges and $j$ bonnets. We can absorb any
part of the graph $G$ that doesn't involve the bonnets using the action of
$\mathcal{O}$.

We only need to consider $\Morph_{\cal{OC}}(o^{\otimes k}, o^{\otimes i}\otimes
t^{\otimes j})$ with $i=0$ and $j=1$, because incoming edges can be
exchanged with outgoing edges and vice versa using the inner product. 
Multiple bonnets must be composites of tori with respect to the open
composition.

%%Embedding $G$ in $\mathbb{R}^3$ so that the bonnet is fixed at the origin
%%and each labelled vertex lies in a distinct plane parallel to the $xy$-plane
%%shows that 
What remains is a composite of open graphs with a single copy of
$B(n)$. 
%%Such an embedding can be obtained by perturbing any embedding that
%%sends the bonnet to 0.
\end{proof}

Let's unwind the definitions in order to determine the chain complex,
$\Torus(A)$, associated to the torus object. Recall that,

$$ (\cal{OC} \otimes_\cal{O} A)(t) = \bigoplus_j \cal{OC}(t,e^{\otimes j}) \otimes A(e^{\otimes j}) = \bigoplus_j \Morph_{\cal{OC}}(e^{\otimes j},t) \otimes A^{\otimes j}$$

modulo the action of $\cal{O}$. This action is determined by the diagram,

$$\begin{diagram}
\cal{OC}(t,e^{\otimes k})\otimes\cal{O}(e^{\otimes j},e^{\otimes k})\otimes A(e^{\otimes j})         & \rTo & \cal{OC}(t,e^{\otimes j})\otimes A(e^{\otimes j}) \\
\dTo        &                   & \dTo\\
\cal{OC}(t,e^{\otimes k})\otimes A(e^{\otimes k}) & \rTo & (\cal{OC}\otimes_\cal{O} A)(t). \\
\end{diagram}$$

As a left $\cal{O}\module$, each $f \in \Morph_{\cal{O}}(e^{\otimes
  j},e^{\otimes k})$ induces a map $f_* : A^{\otimes j} \to A^{\otimes k}$
and, as a right $\cal{O}\module$, each such $f$ induces a map,
$$ f^* : \Morph_{\cal{OC}}(e^{\otimes k}, t) \to \Morph_{\cal{OC}}(e^{\otimes j}, t), $$

given by post-composition. If $g \otimes e^{\otimes k} \in
\Morph_{\cal{OC}}(e^{\otimes k},t)\otimes A^{\otimes k}$ then the diagram above
yields the relation,
$$ f^*(g) \otimes e^{\otimes k} \sim g\otimes f_*(e^{\otimes k}) .$$

Now each complex $\Morph_{\cal{OC}}(e^{\otimes j},t)$ is quasi-isomorphic to
a chain complex of graphs,

$$\Morph_{\cal{OC}}(e^{\otimes j},t) \simeq \bigoplus_g \cal{G}_{(g,j,1)}(\BBar(\Comm)),$$

containing one boundary torus and $j$ edges which, by theorem
\ref{bonnetgen}, is generated by the bonnets $B(n)$.

Recall from section \ref{operadsubsubsec} that two cooperads $\cal{A}$ and
$\cal{B}$ can be isomorphic or quasi-isomorphic. Now the observation that
$\BBar(C) \cong \Lie_{\infty}^* \simeq \Lie^*$, together with lemma
\ref{exactnesslemma} implies that we can think of the complex computing the
relevant homology as graphs with vertices labelled by trees satisfying the
Jacobi (or IHX) relation. Such graphs are $\Comm_{\infty}$ graphs: they
satisfy the shuffle product relation of section \ref{examples of operads} at
each vertex. This can be seen by applying $\Cobar$ to the linear dual of the
short exact sequence $\Comm \to \Ass\to \Lie$. So each equivalence class of
$(\cal{OC}\otimes_{\cal{O}} A)(t)$ under the relation $\sim$ has a unique
representative of the form,
$$\mathbb{Q}\inp{B(n)}\otimes A^{\otimes n}.$$

The differential is determined by the internal differential $\delta$ of
$A$ and the sum of all possible ways to add an edge to the collection of
edges at the vertex of the boundary torus. This can be
described pictorially,
\begin{diagram}
\BPic{bonnet2} &\rTo & \BPic{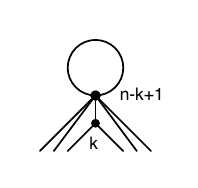}
\end{diagram}

The orientation on the right hand side is induced by the left
hand side. 

In order to describe the object associated to the torus algebraically, we begin by defining $\pTorus(A)$.

$$ \pTorus(A) = \bigoplus_{j=1}^\infty A^{\otimes j} .$$

There is a map $\pi$ from $\pTorus(A)$ onto $(\cal{OC} \otimes_{\cal{O}}
A)(t)$. The kernel of $\pi$ is spanned by shuffles because the
$\Comm_{\infty}$ operad's generators, $m_n$, are precisely those which
vanish on shuffle products. If the shuffle product of tensors is defined by,

$$(a_1 \otimes \cdots \otimes a_i) \ast (a_{i+1} \otimes \cdots \otimes a_n) = \sum_{\sigma \in \Sh(i,n-1)} \pm a_{\sigma(1)} \otimes \cdots \otimes a_{\sigma(n)}, $$

then $\ker(\pi)$ is the ideal of $\pTorus(A)$ generated by the shuffle
products and the object associated to the torus is the chain complex,
$$\Torus(A) = \pTorus(A)/\ker(\pi) .$$

The differential is the sum of the one given by the $\Ass_\infty$ relation,

$$ d (a_1\otimes\cdots\otimes a_n) = \sum_{\substack{i+j = n+1\\i,j \geq 2}} \sum_{s=0}^{n-j} (-1)^{j+s(j+1)} a_1\otimes \cdots\otimes m_j(a_{s+1}\otimes \cdots \otimes a_{s+j+1})\otimes \cdots \otimes a_n$$

and the internal differential of $A$,

$$ \delta(a_1 \otimes \cdots \otimes a_n) = \sum_{i=1}^n a_1\otimes\cdots\otimes \partial (a_i)\otimes \cdots\otimes a_n.$$

If $A$ is a commutative algebra or differential graded commutative algebra,
then the chain complex $\Torus(A)$ agrees with the chain complex computing
Harrison homology, see \cite{MR0220799}. 

The theorem below summarizes the above computation.

\begin{theorem}\label{yaddathm}
If $A$ is a cyclic $\Comm_{\infty}$ algebra and $\mathcal{OC}$ is the
$\Obj(\mathcal{OC})-\mathcal{O}$ bimodule of section \ref{hvwtheoremsec}
then, after identifying $A$ as an $\mathcal{O}$ module, the extension
$\mathcal{OC}\otimes_{\mathcal{O}} A$ associates to the torus object
$t\in\Obj(\mathcal{OC})$ a chain complex, $\Torus(A)$, computing the
Harrison homology of $A$.

$$(\mathcal{OC}\otimes_{\mathcal{O}} A)(t) = \Torus(A) \quad\normaltext{ and }\quad H_*(\Torus(A)) \cong \Har_*(A,A)$$
\end{theorem}

\subsection{Flatness and Exactness}\label{flatexactsec}

In this section we show that the closed category $\mathcal{C}$ acts on the
Harrison complex associated to a torus by theorem \ref{yaddathm}.

For the extension $\mathcal{OC}\otimes_{\mathcal{O}} A$ to be an open-closed
field theory in the sense of definition \ref{octft}, we must show that
$i_*(A)$ is h-split. In order to describe the complex $i_*(A)(t)$, a
simplification can be made,
$$ \cal{OC} \otimes^{\mathbb{L}}_{\cal{O}} A \simeq \cal{OC} \otimes_{\cal{O}} A, $$

by observing that, as an $\Obj(\cal{OC})-\cal{O}$ bimodule, the category
$\cal{OC}$ is flat. 

This is true because there is a natural filtration on the bimodule
$\cal{OC}$ given by the degree of the bonnets. A bonnet with vertex labelled
by $m_n$ must come from a cell of underlying dimension $n-2$. For
instance, the bonnet in degree 0, represented by a trivalent graph, must come
from the trivial forest (or zero dimensional cube), covering only the base
point of the relevant cycle.

Define a filtration $\cal{F}$ of $\cal{OC}$ so that $\cal{F}^0\cal{OC}$
contains the identity elements and the associated graded $\Gr^n \cal{OC}$ is
precisely the $n$th bonnet $B(n)$. Since $d B(n)$ is a sum of bonnets of
lower degree this is a filtration of complexes. There is an induced
filtration on $\cal{OC} \otimes_{\cal{O}} A$ such that the associated graded
$$\Gr^n (\cal{OC} \otimes_{\cal{O}} A)(e^{\otimes i} \otimes t^{\otimes j})$$

consists of placing the identity factors on the $i$ edges and labelling the
$j$ bonnets by elements of $A^{\otimes n}$. Showing that this is true amounts to a
computation nearly identical to that of the previous section.

We will exploit the following lemma,

\begin{lemma}
  If $\varphi : A \to A'$ is a map of filtered complexes such that
  $\varphi_0 : \cal{F}^0 A \to \cal{F}^0 A'$ is a quasi-isomorphism and
  $\varphi_* : \Gr^n A \to \Gr^n A'$ is a quasi-isomorphism then $\varphi_n
  : \cal{F}^n A \to \cal{F}^n A'$ is a quasi-isomorphism for all $n$. In
  particular, $\varphi$ is a quasi-isomorphism.
\end{lemma}

\begin{theorem}
If $A$ is an $h$-split \cal{O} module, then $\mathcal{OC}\otimes_{\mathcal{O}} A$ is an $h$-split $\Obj(\mathcal{OC})$ module.
\end{theorem}

\begin{proof}
We must check that the maps,
$$(\cal{OC} \otimes_{\cal{O}} A)(x) \otimes (\cal{OC} \otimes_{\cal{O}} A)(y) \to (\cal{OC} \otimes_{\cal{O}} A)(x \otimes y)$$

are quasi-isomorphisms. Since this is true in filtration degree 0 it follows
by induction if it holds for the associated graded. A collection of $i$
bonnets labelled by $A$ tensored with a collection of $j$ bonnets labelled
by $A$ is quasi-isomorphic to a collection of $i+j$ bonnets labelled by $A$.
\end{proof}

\begin{theorem}
The category $\mathcal{OC}$ is a flat $\Obj(\cal{OC})-\cal{O}$ bimodule.
That is, the functor $i_* : \cal{O}\module \to \Obj(\cal{OC})\module$ given
by
$$ i_*(A) = \cal{OC} \otimes_{\mathcal{O}} A $$

is exact.
\end{theorem}

\begin{proof}
  Given a quasi-isomorphism of $\Comm_{\infty}$ algebras $\varphi : A \to
  A'$. We must check that the induced map $\cal{OC} \otimes_{\cal{O}} A \to
  \cal{OC} \otimes_{\cal{O}} A'$, is a quasi-isomorphism. Since this is true
  in filtration degree $0$, it follows by induction if it holds for the
  associated graded. The map
$$\bar{\varphi}: \Gr^n(\cal{OC} \otimes_{\cal{O}} A) \to \Gr^n(\cal{OC} \otimes_{\cal{O}} A')$$

is the map between bonnets labelled by tensor powers of $A$ and $A'$. The map $\bar{\varphi}$ is a quasi-isomorphism because a tensor product of quasi-isomorphisms is a quasi-isomorphism.
\end{proof}

\subsection{Deligne's Conjecture}

\begin{corollary}
The category $\cal{C}$ acts on the complex $\Torus(A)$:
$$  \Morph_{\cal{C}}(t^{\otimes i},t^{\otimes j}) \otimes \Torus_*(A)^{\otimes i} \to \Torus_*(A)^{\otimes j}.$$
\end{corollary}
\begin{proof}
  If we consider $A$ as an $\cal{O}\module$ and $\cal{OC}$ as an
  $\cal{OC}-\cal{O}$ bimodule then we can define an $\cal{OC}$ module
  associated to $A$ by $\cal{OC}\otimes_{\cal{O}}^{\mathbb{L}} A$. 
  If $i : \cal{C} \hookrightarrow \cal{OC}$ is the inclusion then
  $i^*(\cal{OC}\otimes_{\cal{O}}A)$ is a $\cal{C}\module$. If $X(A) =
  i^*(\cal{OC}\otimes_{\cal{O}}A)(t)$ is the chain complex associated to the
  torus then there is a natural map
$$ \Morph_{\cal{C}}(t^{\otimes i},t^{\otimes j}) \otimes X(A)^{\otimes i } \to X(A)^{\otimes j}.$$

Earlier, we considered \cal{OC} as an $\Obj(\cal{OC})-\cal{O}$ bimodule and
saw that $\Torus(A) = j^*(\cal{OC}\otimes_{\cal{O}} A)$. On the other hand,
the complex associated to the torus is independent of the choice of
$\Obj(\cal{OC})$ verses $\cal{OC}$. So $X(A)$ is $\Torus(A)$.

\end{proof}

\bibliographystyle{abbrv}
\bibliography{harrison}  %% This looks for the bibliography in jwspinnet.bib

\end{document}